\documentclass[a4paper,11pt,oneside]{article}

\pdfoutput=1

\usepackage{amsmath,amsthm,amsfonts,amssymb,amscd}
\usepackage{mathtools}
\usepackage{float}
\usepackage{subfig} 
\usepackage{enumerate}
\usepackage{booktabs}
\usepackage[noend]{algpseudocode}
\usepackage[normalem]{ulem}
\usepackage{graphicx,bm,xcolor}
\usepackage{algorithm}
\usepackage{algpseudocode}
\usepackage[numbers]{natbib}
\usepackage[hidelinks]{hyperref}
\usepackage{caption}
\usepackage[T1]{fontenc}
\usepackage{t1enc}
\usepackage[polish,german,english]{babel}
\usepackage{verbatim}
\usepackage{diagbox}
\usepackage{colortbl}

\usepackage{url}
\usepackage{authblk}

\algnewcommand{\algorithmicgoto}{\textbf{go to}}
\algnewcommand{\Goto}[1]{\algorithmicgoto~\ref{#1}}
\usepackage[top=2cm, bottom=2cm, left=2cm, right=2cm]{geometry}
\theoremstyle{plain} % default
\newtheorem{theorem}{Theorem}[section]

\newtheorem{remark}[theorem]{Remark}

\newtheorem{form}[theorem]{Formulation}

\theoremstyle{definition} %

\theoremstyle{remark} %

\DeclareMathOperator*{\dG}{dG}
\DeclareMathOperator*{\cG}{cG}
\DeclareMathOperator*{\Ieff}{I_{\text{eff}}}

\newcommand\restrict[1]{\raisebox{-.5ex}{$|$}_{#1}}

\newcommand{\trialc}[1]{{{\textcolor{blue!60!black}{#1}}}}
\newcommand{\testc}[1]{{{\textcolor{green!40!black}{#1}}}}
\newcommand{\puc}[1]{{{\textcolor{brown!40!black}{#1}}}}

\definecolor{brewer1}{HTML}{A6CEE3}
\definecolor{brewer2}{HTML}{1F78B4}
\definecolor{brewer3}{HTML}{B2DF8A}
\definecolor{brewer4}{HTML}{33A02C}

\begin{document}

\title{Tensor-product space-time goal-oriented error control \\ and adaptivity with partition-of-unity dual-weighted residuals  \\ for nonstationary flow problems}
\author[1,2]{Julian Roth}
\author[1]{Jan Philipp Thiele}
\author[3]{Uwe K\"ocher}
\author[1,2]{Thomas Wick}

\affil[1]{Leibniz Universit\"at Hannover, Institut f\"ur Angewandte
  Mathematik, \linebreak AG Wissenschaftliches Rechnen, Welfengarten 1, 30167 Hannover, Germany}
\affil[2]{Universit\'e Paris-Saclay, LMPS - Laboratoire de Mecanique Paris-Saclay,
91190 Gif-sur-Yvette, France}
\affil[3]{Helmut Schmidt University, Faculty of Mechanical Engineering, Holstenhofweg 85, 22043 Hamburg, Germany}

\date{}
\maketitle

%%%%%%%%%%%%%%%%%%%%%%%%%%%%%%%%%%%%%%%%%%%%%%%%%%%%%%%%%%%%
%%                     ABSTRACT                           %%
%%%%%%%%%%%%%%%%%%%%%%%%%%%%%%%%%%%%%%%%%%%%%%%%%%%%%%%%%%%%
\begin{abstract}
\noindent In this work, the dual-weighted residual method is applied to a space-time formulation of nonstationary Stokes and Navier-Stokes flow. 
Tensor-product space-time finite elements are being used to discretize the variational formulation with discontinuous Galerkin finite elements in time and inf-sup stable Taylor-Hood finite element pairs in space. 
To estimate the error in a quantity of interest and drive adaptive refinement in time and space, we demonstrate how the dual-weighted residual method for incompressible flow can be extended to a partition of unity based error localization. 
We derive the space-time Newton method for the Navier-Stokes equations and 
substantiate our methodology on 2D benchmark problems from computational fluid mechanics.
\end{abstract}

%%%%%%%%%%%%%%%%%%%%%%%%%%%%%%%%%%%%%%%%%%%%%%%%%%%%%%%%%%%%
%%                  INTRODUCTION                          %%
%%%%%%%%%%%%%%%%%%%%%%%%%%%%%%%%%%%%%%%%%%%%%%%%%%%%%%%%%%%%
\section{Introduction}
\label{sec_intro}

\noindent This work is devoted to the efficient numerical solution of 
the incompressible Navier-Stokes equations employing 
adaptive Galerkin finite elements in space and time. Adaptivity is based 
on a posteriori goal-oriented error control using 
a partition-of-unity dual-weighted residual (PU-DWR) approach.
The closest studies to ours are on the one hand \cite{BeRa12} as well as the PhD thesis \cite{Schmich2009}
in which similar concepts for full space-time adaptivity for the Navier-Stokes 
equations were developed. The main differences are two-fold that we work
with tensor-product space-time finite elements and that we utilize 
a partition-of-unity for localization (originally proposed for stationary problems in \cite{RiWi15_dwr})
rather than the filtering approach \cite{BraackErn02}.
On the other hand, another closely related study is \cite{ThiWi22_arxiv} in which 
the current PU-DWR method was developed and analyzed for parabolic problems. 

Employing (full) space-time adaptivity using the DWR method dates back to 
\cite{SchmichVexler2008}, who developed and tested the method for 
parabolic problems. Shortly afterwards the space-time DWR method 
was applied to the incompressible Navier-Stokes equations 
in \cite{Schmich2009,BeRa12} and for second-order hyperbolic problems (i.e.,
the elastic wave equation) in \cite{Rade09} and \cite{BaGeRa10}.
In this regard, the first space-time finite element formulation (without adaptivity)
for incompressible flow was proposed in \cite{Tez92}.
Employing space-time (therein called $G2$) for turbulent flow and goal-oriented 
adaptivity for the efficient computation of the mean drag was done in \cite{Hoff09}.

Transport problems with coupled flow with an emphasis on adaptivity 
of the transport part was investigated in \cite{BauBruKoe2021}
and the extension to a multirate framework was suggested in \cite{BruKoeBau2022}.
Goal-oriented error control with space-time formulations, but only applied 
to temporal error control for parabolic problems, the Navier-Stokes 
equations, and fluid-structure interaction was undertaken 
in \cite{MeidnerRichter,MeidnerRichterNavier,FaiWi18}, respectively.
Space-time methods with residual-based error control and adaptivity 
can be found for instance in \cite{LaStein19,Schaf21}.
Finally, \cite{MeiVe07} proposed goal-oriented space-time adaptivity with optimal control 
and \cite{LaSchaf20,LaSchaf21b} developed residual-based adaptivity in optimal control. 

The main objectives of this work are to develop space-time PU-DWR for the 
incompressible Navier-Stokes equations. These developments include 
a detailed computational analysis in form of error reductions and 
studies of effectivity indices. To better understand the behavior, the 
(linear) Stokes equations are considered as well. Therein, 
Galerkin discretizations using tensor products in time and space are employed. 
In time, discontinuous Galerkin $\dG(r)$ elements with polynomial degree $r\geq 0$ 
are used for the primal forward problem. The adjoint problem is approximated 
with globally higher order finite elements. In space the classical Taylor-Hood 
element (i.e., continuous Galerkin $\cG(s)$, with $s\geq 1$ with $s=2$ for the velocities
and $s=1$ for the pressure) 
is employed for the primal problem and a higher-order approximation 
for the adjoint. A delicate issue are pressure-robust discretizations 
on dynamically changing adaptive meshes for which we employ 
a divergence free projection from \cite{BesierWollner2012}, but we make several 
remarks to more modern methods. 
These developments yield 
expressions for error indicators in time and space, which are then utilized 
to formulate a final adaptive algorithm.

The outline of this paper is as follows: In Section \ref{sec_ST}, 
the space-time formulation and discretization are introduced. Next, in 
Section \ref{sec_DWR}, the space-time PU-DWR is developed and error 
estimators and error indicators are derived. Afterward, in Section \ref{sec_tests},
several benchmark tests are adopted to substantiate our algorithmic 
developments. Our work is summarized in Section \ref{sec_conclusions}.

%%%%%%%%%%%%%%%%%%%%%%%%%%%%%%%%%%%%%%%%%%%%%%%%%%%%%%%%%%%%
%%           FORMULATION & DISCRETIZATION                 %%
%%%%%%%%%%%%%%%%%%%%%%%%%%%%%%%%%%%%%%%%%%%%%%%%%%%%%%%%%%%%
\section{Space-time formulation and discretization}
\label{sec_ST}

\noindent We model viscous fluid flow with constant density and temperature by the Navier-Stokes equations \cite{galdi2011, temam2001, Rannacher2000, GiRa1986, turek1999, Glowinski2003}: 
\begin{form}[incompressible isothermal Navier-Stokes equations]
\label{form:navier_stokes}
Find the vector-valued velocity $\bm{v}: (0,T) \times \Omega \rightarrow \mathbb{R}^d$ with $d \in \lbrace 2, 3 \rbrace$ and the scalar-valued pressure $p: (0,T) \times \Omega \rightarrow \mathbb{R}$ such that
\begin{subequations}\label{eq:navier_stokes_strong}
\begin{align}
    \partial_t \bm{v} - \nabla_x \cdot \bm\sigma  + (\bm{v} \cdot \nabla_x)\bm{v}  &= 0 &&\hspace{-1cm} \mathrm{in}\ (0,T) \times \Omega, \\ %
    \nabla_x \cdot \bm{v} &= 0 &&\hspace{-1cm}\mathrm{in}\ (0,T) \times \Omega, \\
    \bm{v} &= \bm{v}_D &&\hspace{-1cm}\mathrm{on}\ (0,T) \times \partial\Omega, \\
    \bm{v}(0) &= \bm{v}^0 &&\hspace{-1cm}\mathrm{in}\ \Omega,
\end{align}
\end{subequations}
where we use the non-symmetric Cauchy stress tensor for the stress $\bm\sigma$ defined as
\begin{align}
    \bm\sigma := \bm\sigma \begin{pmatrix} \bm v \\ p \end{pmatrix} = -pI + \nu \nabla_x \bm{v}.   
\end{align}
Here, $\nu$ is the kinematic viscosity, $\bm{v}^0$ is the initial velocity and $\bm{v}_D$ is a possibly time-dependent Dirichlet boundary value.
\end{form}
\begin{remark}
By omitting the nonlinear convection term $(\bm{v}\cdot \nabla_x)\bm{v}$ we obtain the Stokes problem.
\end{remark}

We employ the following notation for the function spaces. 
Spatial function spaces are denoted by $V$. Without indices $V$ is the analytical function space.
$V_h^{2/1}$ is the space of lowest order Taylor-Hood elements (quadratic FE for velocity, linear FE for pressure) and $V_h^{4/2}$ is the space of fourth order FE for velocity and second order FE for pressure. 
Temporal function spaces are denoted by $X$. Without indices $X$ is the analytical function space.
$X_k^{\dG(r)}$ is the space of discontinuous finite elements in time of $r$.th order and $X_k^{\dG(r+1)}$ is the space of discontinuous finite elements in time of $(r+1)$.th order. 
We denote spatio-temporal function spaces by combining a spatial and a temporal function space:

{\centering\scalebox{0.895}{\parbox{1.115\linewidth}{%
\begin{align}\label{eq:space_time_func_space}
    X(V) := \left\lbrace \begin{pmatrix} \bm v \\ p \end{pmatrix} \, \middle| \, \bm v \in L^2(I, H^1_{0}(\Omega)^d) , \partial_t \bm v \in L^2\left(I, \left(H^1_{0}(\Omega)^d\right)^\ast\right), p \in L^2(I, L^2(\Omega)), \int_\Omega p(t)\ \mathrm{d}\boldsymbol{x} = 0\ \forall t \in I \right\rbrace 
\end{align}
}}}

\noindent will be the space-time Hilbert space on which we define analytical weak formulation. Here $W^* = L(W,\mathbb{R})$ denotes the dual space of $W$, which contains bounded linear operators from $W$ to $\mathbb{R}$. This is a good choice for the analytical function space, since we have a continuous embedding into a function space where the fluid velocity is continuous on $\bar{I}$ \cite{BeRa12, dautray1999}. This ensures that the initial condition is well-defined for the velocity.
\begin{remark}
 When we have homogeneous Neumann boundary conditions $\bm \sigma \cdot \bm n = \bm 0$ on $\Gamma_{\text{out}} \subset \partial \Omega$, the so-called do-nothing (outflow) condition \cite{Heywood_Rannacher_Turek_1996},
the pressure does not require normalization for its uniqueness and we have the spatio-temporal function space
\begin{align}
    X(V) := \left\lbrace \begin{pmatrix} \bm v \\ p \end{pmatrix} \, \middle| \, \bm v \in L^2(I, H^1_{0, \Gamma_D}(\Omega)^d) , \partial_t \bm v \in L^2\left(I, \left(H^1_{0, \Gamma_D}(\Omega)^d\right)^\ast\right), p \in L^2(I, L^2(\Omega)) \right\rbrace,  
\end{align}
where $\Gamma_D = \partial \Omega \setminus \Gamma_{\text{out}}$.
\end{remark}

Finally, to simplify notation we introduce
\begin{align}
    (f,g) := (f,g)_{L^2(\Omega)} := \int_\Omega f \cdot g\ \mathrm{d}\bm x, \qquad (\!(f,g)\!) := (f,g)_{L^2(I, L^2(\Omega))} := \int_I (f, g)\ \mathrm{d}t.
\end{align}
In this notation, $f \cdot g$ represents the Euclidean inner product if $f$ and $g$ are scalar- or vector-valued and the Frobenius inner product if $f$ and $g$ are matrices.

For a given initial value $\bm v^0 \in L^2(\Omega)^d$ the weak formulation reads: Find $\bm U := \begin{pmatrix} \bm v \\ p \end{pmatrix} \in \begin{pmatrix} \bm v_D \\ 0 \end{pmatrix} + X(V)$ such that
\begin{align}
    A(\bm U)(\bm \Phi) := (\!(\partial_t \bm v, \bm \Phi^v)\!) + a(\bm U)(\bm \Phi) + (\bm v(0) - \bm v^0, \bm \Phi^v(0)) = 0 \qquad \forall \bm \Phi := \begin{pmatrix} \bm \Phi^v \\ \Phi^p \end{pmatrix} \in X(V)
\end{align}
where for the Navier-Stokes equations, we have
\begin{align}
    a(\bm U)(\bm \Phi) :=  - (\!(p, \nabla_x \cdot \bm \Phi^v)\!) + \nu (\!(\nabla_x \bm v, \nabla_x \bm \Phi^v)\!) + (\!((\bm{v} \cdot \nabla_x)\bm{v}, \bm \Phi^v)\!) + (\!(\nabla_x \cdot \bm v, \Phi^p)\!),
\end{align}
and for the Stokes equations, we have
\begin{align}
    a(\bm U)(\bm \Phi) :=  - (\!(p, \nabla_x \cdot \bm \Phi^v)\!) + \nu (\!(\nabla_x \bm v, \nabla_x \bm \Phi^v)\!) + (\!(\nabla_x \cdot \bm v, \Phi^p)\!).
\end{align}

\subsection{Discretization in time}
\label{sec:time_discretization}
\noindent Let $ \bar{I} = [0,T] = \lbrace 0 \rbrace \cup \bigcup_{m = 1}^M I_m $,
with $I_m := (t_{m-1},t_m]$ be a partitioning of time. Then, we define the semi-discrete space as

\begin{figure}[H]
    \centering
    \subfloat[dG(0)]{
        \includegraphics{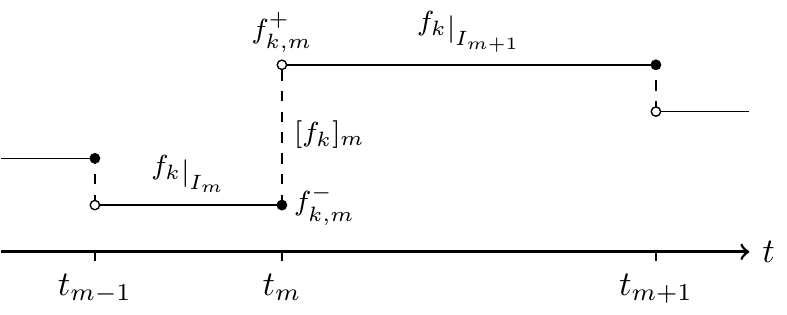}
    }%
    \quad
    \subfloat[dG(1)]{
        \includegraphics{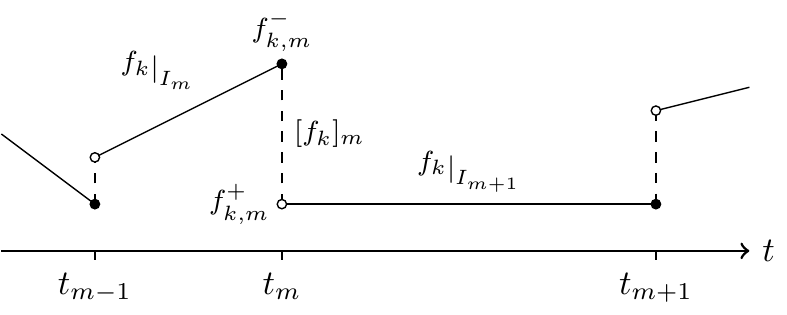}
    }%
    \caption{dG(r) time discretization}
    \label{fig:dGr_time_discretization}
\end{figure}

{\centering\scalebox{0.75}{\parbox{1.325\linewidth}{%
\begin{align}
    X_k^{\dG(r)}(V) := \left\lbrace \begin{pmatrix} \bm{v_k} \\ p_k \end{pmatrix} \middle| \bm{v_k}(0) \in L^2(\Omega)^d, \bm{v_k}\restrict{I_m} \in P_r(I_m, H^1_{0}(\Omega)^d), p_k\restrict{I_m} \in P_r(I_m, L^2(\Omega))\ \forall 1 \leq m \leq M, \int_\Omega p_k(t)\ \mathrm{d}\boldsymbol{x} = 0\ \forall t \in I \right\rbrace,
\end{align}
}}}

\noindent where the space-time function space (\ref{eq:space_time_func_space}) has been discretized in time with the discontinuous Galerkin method of order $r\in \mathbb{N}_0$ ($\dG(r)$). $P_r(I_m, Y)$ is the space of polynomials of order $r$, which map from the time interval $I_m$ into the space $Y$. Since functions in $X_k^{\dG(r)}(V)$ can have discontinuities between the time intervals, we define the limits of $f_k$ at time $t_m$ from above and from below for a function $f_k$
\begin{align}
    f_{k,m}^\pm := \lim_{\epsilon \searrow 0} f_k(t_m \pm \epsilon),
\end{align}
and the jump of the function value of $f_k$ at time $t_m$ as
\begin{align}
    [f_k]_m := f_{k,m}^+ - f_{k,m}^-.
\end{align}
The $\dG(r)$ time discretization for the cases $r = 0$ and $r = 1$ is illustrated in Figure \ref{fig:dGr_time_discretization}.

\noindent It has been derived in \cite{Schmich2009} for the Navier-Stokes equations that we obtain the backward Euler scheme by using suitable numerical quadrature for the temporal integrals of the $\dG(0)$ discretization.
Therefore, it can be seen as a variant of that scheme. 
Additionally, the corresponding time-stepping scheme of $\dG(1)$ has been formulated therein. The $\dG(r)$ time discretization of the Navier-Stokes equations reads: \\ 
\indent Find $\bm{U_k} := \begin{pmatrix} \bm{v_k} \\ p_k \end{pmatrix} \in \begin{pmatrix} \bm{v}_D \\ 0 \end{pmatrix} + X_k^{\dG(r)}(V)$ such that
\begin{equation}
    \begin{aligned}
    A_k(\bm{U_k})(\bm{\Phi_k}) := &\sum_{m = 1}^M \int_{I_m}(\partial_t \bm{v_k}, \bm{\Phi_k}^v) - (p_k, \nabla_x \cdot \bm{\Phi_k}^v) + \nu (\nabla_x \bm{v_k}, \nabla_x \bm{\Phi_k}^v)\ \mathrm{d}t \\
    &+ \sum_{m = 1}^M \int_{I_m} ((\bm{v_{\bm{k}}} \cdot \nabla_x)\bm{ v_{\bm{k}}}, \bm \Phi^v_{\bm{k}}) + (\nabla_x \cdot \bm{v_{\bm{k}}}, \Phi_{k}^p)\ \mathrm{d}t \\
    &+ \sum_{m=0}^{M-1} ([\bm{v_k}]_m, \bm{\Phi}_{\bm{k},m}^{v,+})  + (\bm{v}_{\bm{k},0}^-(0) - \bm v^0, \bm{\Phi}_{\bm{k},0}^{v,-}(0)) = 0 \\
    &\forall \bm{\Phi_k} := \begin{pmatrix} \bm{\Phi_k}^v \\ \Phi_k^p \end{pmatrix} \in X_k^{\dG(r)}(V).
    \end{aligned}   
\end{equation}

\subsection{Discretization in space}
Next, we describe the spatial discretization of the $\dG(r)$ formulation for the (Navier-)Stokes equations with continuous finite elements. Let $\mathbb{T}_h$ be a shape-regular mesh of the  domain $\Omega \subset \mathbb{R}^d$. The elements $K \in \mathbb{T}_h$ are a non-overlapping covering of $\bar{\Omega}$, where the elements $K$ are quadrilaterals for $d = 2$ and hexahedrals for $d = 3$. The discretization parameter $h$ is defined element-wise as $h\restrict{K} = h_K = \operatorname{diam}(K)$. After adaptive refinement in space, $h_K$ can vary across elements and we need to account for hanging nodes in the assembly of the matrices and vectors. These are degrees of freedom that are located in the middle of the neighboring elements' edge ($d=2$) or face ($d=3$). Furthermore, on different time intervals we allow for changing spatial meshes, i.e. possibly $\mathbb{T}_h(I_m) \neq \mathbb{T}_h(I_n)$, but for simplicity we drop the dependence on the time interval and only mention it when it is not obvious from the context. 

During the spatial discretization of the mixed system of the (Navier-)Stokes equations, we need to account for the inf-sup or  Ladyzhenskaya-Babuska-Brezzi (LBB) condition
\begin{align}
    \underset{q \in L^2(\Omega)}{\inf}\ \underset{\bm\varphi \in H^1_{0}(\Omega)^d}{\sup}\ \frac{(q, \nabla \cdot \bm\varphi)}{\|q\|_{L^2(\Omega)} \| \bm\varphi\|_{H^1(\Omega)^d}} \geq \gamma > 0,
\end{align}
which guarantees well-posedness (existence, uniqueness and stability) of the 
pressure solution under suitable boundary conditions and assumptions
on the domain regularity; see e.g., \cite{GiRa1986,Ra17_3}.
In the spatial discretization we now replace $L^2(\Omega)$ and $H^1_{0}(\Omega)^d$ with conforming finite element spaces that still fulfill a discrete version of this stability estimate. For this we define the mixed finite element space
\begin{align}
    V_h^{s_v/s_p} :=  V_h^{s_v/s_p}(\mathbb{T}_h) := \left\lbrace \begin{pmatrix} \bm v \\ p \end{pmatrix} \in  C(\bar \Omega)^{d+1} \middle| \bm{v}\restrict{K} \in \mathcal{Q}_{s_v}(K)^d, p\restrict{K} \in \mathcal{Q}_{s_p}(K) \quad  \forall K \in \mathbb{T}_h \right\rbrace,
\end{align}
where $\mathcal{Q}_{s}(K)$ is being constructed by mapping tensor-product polynomials of degree $s$ from the master element $(0,1)^d$ to the element $K$. In our computations, we use the function space $V_h^{2/1}$, the so-called Taylor-Hood elements \cite{TaylorHood1974}  with quadratic finite elements for the velocity and linear finite elements for the pressure. Additionally, we use as $V_h^{4/2}$ an enriched space with fourth order finite elements for the velocity and second order finite elements for the pressure.
In $\mathbb{R}^2$, we could have also used $V_h^{3/2}$ as an enriched space, since $\mathcal{Q}_{3}/\mathcal{Q}_{2}$ was shown to be an inf-sup stable Taylor-Hood pair \cite{Stenberg1990}, but this must not be the case in 3D anymore.
Hence, we approximate the enriched velocity by $\mathcal{Q}_{4}$ instead of  $\mathcal{Q}_{3}$ finite elements in space.
The fully discretized problem then reads: Find $\bm U_{\bm{kh}} := \begin{pmatrix} \bm{v_{\bm{kh}}} \\ p_{kh} \end{pmatrix} \in \begin{pmatrix} \bm{v}_D \\ 0 \end{pmatrix} + X_k^{\dG(r)}\left(V_h^{s_v/s_p}\right)$ such that
\begin{equation}
\begin{aligned}
\label{eq:space_time_discrete_formulation}
    A_{kh}(\bm U_{\bm{kh}})(\bm{\Phi_{\bm{kh}}} ) &:= \sum_{m = 1}^M \int_{I_m}(\partial_t \bm{v_{\bm{kh}}}, \bm{\Phi_{\bm{kh}}}^v) - (p_{kh}, \nabla_x \cdot \bm{\Phi_{\bm{kh}}}^v) + \nu (\nabla_x \bm{v_{\bm{kh}}}, \nabla_x \bm{\Phi_{\bm{kh}}}^v)\ \mathrm{d}t   \\
    &\hspace{-2cm}+ \sum_{m = 1}^M \int_{I_m} ((\bm{v_{\bm{kh}}} \cdot \nabla_x)\bm{ v_{\bm{kh}}}, \bm \Phi^v_{\bm{kh}}) + (\nabla_x \cdot \bm{v_{\bm{kh}}}, \Phi_{kh}^p)\ \mathrm{d}t \\
    &\hspace{-2cm}+ \sum_{m=0}^{M-1} ([\bm{v_{\bm{kh}}}]_m, \bm{\Phi}_{\bm{kh},m}^{v,+})  + (\bm{v}_{\bm{kh},0}^-(0) - \bm v^0, \bm{\Phi}_{\bm{kh},0}^{v,-}(0)) = 0 \\
    &\forall \bm{\Phi_{\bm{kh}}} := \begin{pmatrix} \bm{\Phi_{\bm{kh}}}^v \\ \Phi_{kh}^p \end{pmatrix} \in X_k^{\dG(r)}\left(V_h^{s_v/s_p}\right). 
\end{aligned}
\end{equation}

\subsection{Tensor-product space-time discretization}
\label{sec:tensor_product_space_time}
The above space-time analysis is required for the theory of error estimation for the time-dependent partial differential equations.
However, this formulation is usually transformed into a time-stepping scheme in the literature by applying a suitable numerical quadrature to the temporal integrals (see also Section \ref{sec:time_discretization}). 
Instead we will follow the approach of Bruchh\"auser et al. \cite{BruKoeBau2022}[Sec. 4], which uses tensor-product spaces to directly work with the fully discrete weak formulation (\ref{eq:space_time_discrete_formulation}). In this approach the space-time cylinder $(0,T) \times \Omega$ is divided into space-time slabs
\begin{align}
    Q_n := \left( \bigcup\limits_{m = m_n}^{m_{n+1}-1} I_m \right) \times \Omega, \qquad 1 = m_1 < m_2 < ... < m_l < m_{l+1}-1 = M,    
\end{align}
on which the weak formulation is solved sequentially.
The spatial mesh is slabwise constant, i.e. $\mathbb{T}_h(I_{m_n}) =  \dots = \mathbb{T}_h(I_{m_{n+1}-1})$. 
Hence, we can define the space-time tensor-product finite element basis on $Q_n$ by taking the tensor-product of the temporal discontinuous Galerkin basis functions and the spatial continuous Galerkin basis functions. 
This allows for the extension of a simple finite element code for a stationary Stokes problem by creating a temporal mesh and adding a loop over the temporal degrees of freedom in the assembly of the FEM matrices. 
Furthermore, unlike for the time-stepping formulations of (\ref{eq:space_time_discrete_formulation}) where each polynomial degree in time needs to be derived and implemented separately, changing the temporal degree of the space-time discretization can be performed simply by changing the polynomial degree of the temporal finite elements and our code theoretically allows for hp-adaptivity in time, since the temporal degree could be chosen on each slab individually.
Another benefit of using space-time slabs for adaptive refinement is that the spatial mesh is fixed for the length of the slab which is sometimes desirable for numerical stability.
On the other hand, using slabs with multiple time elements comes at the cost of solving larger linear equation systems.

\subsection{Linearization of the semilinear form}
\label{sec:linearization_nse}
For the Navier-Stokes equations, we need to include the nonlinear convection term $(\!((\bm{v} \cdot \nabla_x)\bm{v}, \bm \Phi^v)\!)$ in our variational formulation $A_{kh}(\bm U_{\bm{kh}})(\bm{\Phi_{\bm{kh}}}) = 0$. 
Therefore, we need to apply linearization to obtain a solvable linear equation system.
In this paper, we decided to solve this variational root finding problem by using the residual-based Newton's method.
The resulting nonlinear iteration step $j+1\in\mathbb{N}$ reads as:\\
Solve the problem:
\\
\indent Find $\bm{\delta U}_{\bm{kh}} := \begin{pmatrix} \bm{\delta v_{\bm{kh}}} \\ \delta p_{kh} \end{pmatrix} \in X_k^{\dG(r)}\left(V_h^{s_v/s_p}\right)$ such that
\begin{align}
    A^\prime_{kh,\bm{U}}(\bm{U^j_{\bm{kh}}})(\bm{\delta U_{\bm{kh}}}, \bm{\Phi_{\bm{kh}}})  = -A_{kh}(\bm U^j_{\bm{kh}})(\bm{\Phi_{\bm{kh}}}) \qquad \forall \bm{\Phi_{\bm{kh}}} \in X_k^{\dG(r)}\left(V_h^{s_v/s_p}\right).
\end{align}
Obtain the new iterate by:
\begin{align}
    \bm{U^{j+1}_{\bm{kh}}} = \bm{U^j_{\bm{kh}}}+\bm{\delta U_{\bm{kh}}}
\end{align}
\begin{remark}
    Note that we will also need to assemble the directional derivative $A^\prime_{kh,\bm{U}}$ for the adjoint problem (\ref{eq:space_time_discrete_general_adjoint}) with the minor difference that there the trial and test functions are interchanged. This corresponds to transposing the system matrix from a Newton step.
\end{remark}
Hence, we now need to compute the Fréchet derivative of the semi-linear form $A_{kh}(\bm U_{\bm{kh}})(\bm{\Phi_{\bm{kh}}})$ with respect to $\bm{U}$ in the direction of $\bm{\delta U_{\bm{kh}}}$. For a linear functional $f: X_k^{\dG(r)}\left(V_h^{s_v/s_p}\right) \rightarrow \mathbb{R}$ 
\begin{align}
    f^\prime(\bm{U})(\bm{\delta U}) = f(\bm{\delta U}),
\end{align}
holds.
Hence, it remains to compute the Fréchet derivative of the convection term
\begin{align}
    f: X_k^{\dG(r)}\left(V_h^{s_v/s_p}\right) \rightarrow \mathbb{R},\quad \bm{U} := \begin{pmatrix} \bm v \\ p \end{pmatrix} \mapsto (\!((\bm{v} \cdot \nabla_x)\bm{v}, \bm \Phi^v)\!),
\end{align}
which directly follows from the linearity of differentiation as
\begin{align}
    f^\prime(\bm{U})(\bm{\delta U}) = (\!((\bm{\delta v} \cdot \nabla_x)\bm{v} + (\bm{v} \cdot \nabla_x)\bm{\delta v}, \bm \Phi^v)\!).
\end{align}
The derivative of the complete semi-linear form is thus given by
\begin{equation}
\begin{aligned}
    &A^\prime_{kh,\bm{U}}(\bm{U_{\bm{kh}}})(\bm{\delta U_{\bm{kh}}}, \bm{\Phi_{\bm{kh}}})  = \sum_{m = 1}^M \int_{I_m}(\partial_t \bm{\delta v_{\bm{kh}}}, \bm{\Phi_{\bm{kh}}}^v) - (\delta p_{kh}, \nabla_x \cdot \bm{\Phi_{\bm{kh}}}^v) + \nu (\nabla_x \bm{\delta v_{\bm{kh}}}, \nabla_x \bm{\Phi_{\bm{kh}}}^v) \ \mathrm{d}t \\
    &\hspace{1cm}+ \sum_{m = 1}^M \int_{I_m} ((\bm{\delta v_{\bm{kh}}} \cdot \nabla_x)\bm{v_{\bm{kh}}} + (\bm{v_{\bm{kh}}} \cdot \nabla_x)\bm{\delta v_{\bm{kh}}}, \bm \Phi^v_{\bm{kh}}) + (\nabla_x \cdot \bm{\delta v_{\bm{kh}}}, \Phi_{kh}^p)\ \mathrm{d}t \\
    &\hspace{1cm}+ \sum_{m=0}^{M-1} ([\bm{\delta v_{\bm{kh}}}]_m, \bm{\Phi}_{\bm{kh},m}^{v,+})  + (\bm{\delta v}_{\bm{kh},0}^-(0) - \bm{\delta v}^0, \bm{\Phi}_{\bm{kh},0}^{v,-}(0)).
\end{aligned}
\end{equation}
Finally, we can formulate the full algorithm of Newton's method, which we use for the numerical experiments.  

\begin{algorithm}[H]
\caption{Slabwise residual-based Newton method with backtracking line search} \label{algo:newton_method}
\hspace*{\algorithmicindent} \textbf{Parameters:} 
$TOL\in\mathbb{R}_{+}$, $\lambda\in(0,1]$, $\max_\text{line search}\in\mathbb{N}_{+}$, $\max_\text{Newton}\in\mathbb{N}_{+}$\\
\hspace*{\algorithmicindent} \textbf{Input:} Initial guess $\bm{U_{\bm{kh}}}^{m,0} \in \begin{pmatrix} \bm{v}_D \\ 0 \end{pmatrix} + X_k^{\dG(r)}\left(V_h^{s_v/s_p}\right)$ on slab $Q_m$ \\
\hspace*{\algorithmicindent} \textbf{Output:} FEM solution $\bm{U_{\bm{kh}}}^{m} \in \begin{pmatrix} \bm{v}_D \\ 0 \end{pmatrix} + X_k^{\dG(r)}\left(V_h^{s_v/s_p}\right)$ on slab $Q_m$.
\begin{algorithmic}[1]
\For{$j = 0, 1, 2, ...,\max_\text{Newton}$}
    \State Find $\bm{\delta U_{\bm{kh}}}^{m} \in X_k^{\dG(r)}\left(V_h^{s_v/s_p}\right)$ such that \label{line:newton2}
    \begin{align}
        A^\prime_{kh,\bm{U}}(\bm{U_{\bm{kh}}}^{m,j})(\bm{\delta U_{\bm{kh}}}^m, \bm{\Phi_{\bm{kh}}})  = -A_{kh}(\bm U_{\bm{kh}}^{m,j})(\bm{\Phi_{\bm{kh}}}) \qquad \forall \bm{\Phi_{\bm{kh}}} \in X_k^{\dG(r)}\left(V_h^{s_v/s_p}\right)
    \end{align}
    \For{$k=0,1,\dots,\max_\text{line search}$} \Comment{line search}
    \State Compute new solution candidate
    \begin{align}
        \bm U_{\bm{kh}}^{m,j+1} &= \bm U_{\bm{kh}}^{m,j} + \lambda^k\ \bm{\delta U_{\bm{kh}}}^{m}
    \end{align}
    \State Check line search criterion
    \begin{align}\label{eq:convergence_criterion2}
        \| R(\bm U_{\bm{kh}}^{m,j+1}) \| \leq \| R(\bm U_{\bm{kh}}^{m,j}) \|
    \end{align}
    \If{(\ref{eq:convergence_criterion2}) is fullfilled}
        \State stop line search
    \EndIf
    \EndFor
    \State Check the stopping criterion
        \begin{align}\label{eq:stopping_criterion2}
            \| R(\bm U_{\bm{kh}}^{m,j+1}) \| \leq TOL.
        \end{align}
    \If{(\ref{eq:stopping_criterion2}) is fulfilled}
        \State return $\bm U_{\bm{kh}}^{m} := \bm U_{\bm{kh}}^{m,j+1}$.
    \EndIf
\EndFor
\end{algorithmic}
\end{algorithm}
\noindent For the numerical tests, we used $\max_{\text{line search}} = 10$ as the maximum number of line search steps and $\lambda = 0.6$ as the line search damping parameter. We set the maximum number of Newton steps to $\max_{\text{Newton}}=10$ and the tolerance $TOL$ to $10^{-10}$.

Since we are dealing with a time-dependent problem, the initial guess can be chosen as the solution of the previous slab evaluated at its end time, i.e. for the slabs $Q_{m-1} := (t_{m-2},t_{m-1}] \times \Omega$ and $Q_{m} := (t_{m-1},t_{m}] \times \Omega$, we define
\begin{align}
    \bm{U_{\bm{kh}}}^{m,0} \equiv \bm{U_{\bm{kh}}}^{m-1}(t_{m-1}). 
\end{align}
\begin{remark}
    Alternatively one could  create the initial guess $\bm{U_{\bm{kh}}}^{m,0}$ as the temporal 
    extrapolant of the entire solution from the previous slab, e.g. $\bm{U_{\bm{kh}}}^{m,0}(t) := \bm{U_{\bm{kh}}}^{m-1}(t_{m-2}) + (t - t_{m-2}) \cdot \dfrac{\bm{U_{\bm{kh}}}^{m-1}(t_{m-1})-\bm{U_{\bm{kh}}}^{m-1}(t_{m-2})}{t_{m-1}-t_{m-2}}$ for linear extrapolation in time.
\end{remark}
Additionally, we need to prescribe the non-homogeneous Dirichlet conditions on the initial guess. However, the boundary conditions for the update $\bm{\delta U_{\bm{kh}}}^m$ become homogeneous, i.e. we enforce $\bm{\delta U_{\bm{kh}}}^m = \bm{0}$ on the Dirichlet boundary.

We observe that in Algorithm \ref{algo:newton_method} we need to compute the directional derivative $A^\prime_{kh,\bm{U}}(\bm{U_{\bm{kh}}}^{m,j})(\bm{\delta U_{\bm{kh}}}^m, \bm{\Phi_{\bm{kh}}})$ in each step of the Newton method. This is computationally expensive, since we need to assemble the Jacobian matrix in every step. Instead, in our program we use the approximation
\begin{align}
    A^\prime_{kh,\bm{U}}(\bm{U_{\bm{kh}}}^{m,j+1})(\cdot, \bm{\Phi_{\bm{kh}}}) \approx A^\prime_{kh,\bm{U}}(\bm{U_{\bm{kh}}}^{m,j})(\cdot, \bm{\Phi_{\bm{kh}}}),
\end{align}
when the reduction rate
\begin{align}
    \theta_j := \frac{\| R(\bm U_{\bm{kh}}^{m,j+1}) \|}{\| R(\bm U_{\bm{kh}}^{m,j}) \|}
\end{align}
is sufficiently small. We use these so called simplified-Newton steps, when $\theta_j \leq \theta_{\text{max}} := 0.1$.
\subsection{Divergence free projection}
\label{sec:div_free_projection}

\noindent In fluid dynamics, working with dynamically changing meshes leads to non-physical oscillations in the values of goal functionals, which has been proven in \cite{BesierWollner2012} for a Stokes model problem. Further analysis of the numerical artifacts caused by the violation of the inf-sup condition, due to the usage of time varying spatial meshes can be found in \cite{BraackLangTaschenberger2013, BaenschKarakatsaniMakridakis2018, AlauzetMehrenberger2010, LedererLinke2017}. To avoid the resulting numerical artifacts, Besier and Wollner \cite{BesierWollner2012} proposed to either repeat the last time-step on the finite element space of the new slab or to use a divergence free projection for the initial condition when switching to a different spatial mesh. Here, we decided to focus on the latter and computationally analyze the impact of a divergence free $L^2$ and a divergence free $H^1_0$ projection \cite{BesierWollner2012}. 

\begin{remark}
    We notice that this approach uses a global correction, which is expensive and requires the solution of an additional problem. Using the Helmholtz-decomposition, Linke \cite{linke2014role} showed that the continuous 
    velocity solution is invariant to irrotational forcings. However, conforming mixed finite elements 
    lose these properties on the discrete level leading to non-physical solutions.
    Various fixes were proposed and studied using local projections to $H(div)$-conforming, i.e. discretely divergence free, finite elements inside the weak formulation of the (Navier-)Stokes equations \cite{linke2014role},\cite{linke2016robust},\cite{doi:10.1137/15M1047696},\cite{LedererLinke2017}.
    However, these studies for fluid flow have all been done on triangular elements and applying them with corresponding
    quadrilateral elements did not yield the desired results
    for our fluid flow formulations. In incompressible solid mechanics, quadrilateral elements with such pressure-robust schemes seem to work though \cite{BaMaWaWiWo22,BaWo22}.
\end{remark}

\noindent This problem can be formalized as follows: \\
\indent Let $Q_m := (t_{m-1}, t_m] \times \Omega$ and $Q_{m+1} := (t_{m}, t_{m+1}] \times \Omega$ be two space-time slabs with different spatial triangulations $\mathbb{T}_h^m := \mathbb{T}_h(Q_m) \neq \mathbb{T}_h(Q_{m+1}) =: \mathbb{T}_h^{m+1}$. Then, we need to find a projection $\Pi$ of the initial condition $\bm U_{\bm{kh}}(t_m) := \begin{pmatrix} \bm{v_{\bm{kh}}}(t_m) \\ p_{kh}(t_m) \end{pmatrix} \in \begin{pmatrix} \bm{v}_D(t_m) \\ 0 \end{pmatrix} + V_h^{s_v/s_p}(\mathbb{T}_h^m) $ 
into the function space of the next slab 
$\begin{pmatrix} \bm{v}_D(t_m) \\ 0 \end{pmatrix} + V_h^{s_v/s_p}(\mathbb{T}_h^{m+1})$ such that the projected initial condition is divergence free, i.e. 
\begin{align}
    (\nabla_x \cdot (\Pi\bm{v_{\bm{kh}}}(t_m)), \Phi_{h}^p) = 0 \qquad \forall \Phi_{h}^p \in  V_h^{s_p}(\mathbb{T}_h^{m+1}).
\end{align}
To achieve this goal Besier and Wollner \cite{BesierWollner2012} presented and analyzed a divergence free $L^2$ projection and a divergence free $H^1_0$ projection. \\

\noindent \textbf{Divergence free $L^2$ projection} \\
Find $\Pi\bm U_{\bm{kh}}(t_m) = \begin{pmatrix} \Pi\bm{v_{\bm{kh}}}(t_m) \\ \Pi p_{kh}(t_m) \end{pmatrix} \in \begin{pmatrix} \bm{v}_D(t_m) \\ 0 \end{pmatrix} + V_h^{s_v/s_p}(\mathbb{T}_h^{m+1})$ such that
\begin{equation}
\begin{aligned}
    (\Pi\bm v_{\bm{kh}}(t_m), \bm{\Phi_{h}}^v) - (\Pi p_{\bm{kh}}(t_m), \nabla_x \cdot \bm{\Phi_{h}}^v) + (\nabla_x \cdot (\Pi\bm v_{\bm{kh}}(t_m)), \Phi_{h}^p) =& (\bm v_{\bm{kh}}(t_m), \bm{\Phi_{h}}^v) \hspace{2cm}\\  \forall \bm{\Phi_{\bm{h}}} :=& \begin{pmatrix} \bm{\Phi_{\bm{h}}}^v \\ \Phi_{h}^p \end{pmatrix} \in V_h^{s_v/s_p}(\mathbb{T}_h^{m+1}).
\end{aligned}
\end{equation}

\noindent \textbf{Divergence free $H^1_0$ projection} \\
Find $\Pi\bm U_{\bm{kh}}(t_m) = \begin{pmatrix} \Pi\bm{v_{\bm{kh}}}(t_m) \\ \Pi p_{kh}(t_m) \end{pmatrix} \in \begin{pmatrix} \bm{v}_D(t_m) \\ 0 \end{pmatrix} + V_h^{s_v/s_p}(\mathbb{T}_h^{m+1})$ such that
\begin{equation}
\footnotesize
\begin{aligned}
    (\nabla_x (\Pi\bm v_{\bm{kh}}(t_m)), \nabla_x \bm{\Phi_{h}}^v) - (\Pi p_{\bm{kh}}(t_m), \nabla_x \cdot \bm{\Phi_{h}}^v) + (\nabla_x \cdot (\Pi\bm v_{\bm{kh}}(t_m)), \Phi_{h}^p) =&
     (\nabla_x \bm v_{\bm{kh}}(t_m), \nabla_x \bm{\Phi_{h}}^v) \hspace{2cm}\\
\forall \bm{\Phi_{\bm{h}}} :=& \begin{pmatrix} \bm{\Phi_{\bm{h}}}^v \\ \Phi_{h}^p \end{pmatrix} \in 
V_h^{s_v/s_p}(\mathbb{T}_h^{m+1}).
\end{aligned}
\end{equation}

%%%%%%%%%%%%%%%%%%%%%%%%%%%%%%%%%%%%%%%%%%%%%%%%%%%%%%%%%%%%
%%                SPACE-TIME DWR                          %%
%%%%%%%%%%%%%%%%%%%%%%%%%%%%%%%%%%%%%%%%%%%%%%%%%%%%%%%%%%%%
\section{Space-time dual weighted residual error estimation}
\label{sec_DWR}

In many practical applications, we are in general not interested in the entire solution of the Navier-Stokes equations but only in some quantity of interest. In this chapter, we will describe how we can use the Dual Weighted Residual (DWR) method 
\cite{becker_rannacher_2001,BeRa96} (see also other references
\cite{bangerth_rannacher_2003, SchmichVexler2008, BeRa12,AinsworthOden:2000,Od18})
to estimate the overall error of our finite element solution in this quantity of interest and how we can localize the error with a space-time partition of unity (PU) \cite{RiWi15_dwr, ThiWi22_arxiv} to refine the spatial and temporal meshes.

\subsection{DWR for parabolic PDEs}
\label{sec:DWR_parabolic}
In the following analysis, we will work with homogeneous boundary conditions to simplify the presentation.
Let us consider a goal functional $J: X(V) + X_k^{\dG(r)}(V) \rightarrow \mathbb{R}$ of the form
\begin{align}
    J(\bm{U}) = \int_0^T J_1(\bm{U}(t))\ \mathrm{d}t + J_2(\bm{U}(T)), 
\end{align}
which represents our physical quantity of interest.
We now want to minimize the difference between the quantity of interest of the analytical solution $\bm{U}$ and its  
FEM approximation $\bm{U}_{\bm{kh}}$, i.e.
\begin{align}
    J(\bm{U}) - J(\bm{U}_{\bm{kh}})
\end{align}
subject to the constraint that the variational formulation is being satisfied. 
As in constrained optimization, we define the Lagrange functionals
\begin{align}
    &\mathcal{L}: X(V) \times X(V) \rightarrow \mathbb{R}, &&(\bm{U}, \bm{Z}) \mapsto J(\bm{U}) - A(\bm{U})(\bm{Z}),\\ 
    &\mathcal{L}_k: X_k^{\dG(r)}(V) \times X_k^{\dG(r)}(V) \rightarrow \mathbb{R},
    &&(\bm{U_k}, \bm{Z_k}) \mapsto J(\bm{U_k}) - A_k(\bm{U_k})(\bm{Z_k}), \\ 
    &\mathcal{L}_{kh}: X_k^{\dG(r)}(V_h^{s_v/s_p}) \times X_k^{\dG(r)}(V_h^{s_v/s_p}) \rightarrow \mathbb{R},
    &&(\bm{U_{\bm{kh}}}, \bm{Z_{\bm{kh}}}) \mapsto J(\bm{U_{\bm{kh}}}) - A_{kh}(\bm{U_{\bm{kh}}})(\bm{Z_{\bm{kh}}}). 
\end{align}
The stationary points $(\bm{U}, \bm{Z})$, $(\bm{U_k}, \bm{Z_k})$ and $(\bm{U}_{\bm{kh}}, \bm{Z}_{\bm{kh}})$ of the Lagrange functionals $\mathcal{L}$, $\mathcal{L}_k$ and $\mathcal{L}_{kh}$ need to satisfy the Karush-Kuhn-Tucker first order optimality conditions. Firstly, the stationary points are solutions to the equations
\begin{align}
    \mathcal{L}^\prime_{\bm{Z}}(\bm{U}, \bm{Z})(\bm{\delta Z}) &= 0 \quad \forall \bm{\delta  Z} \in X(V), \\
    \mathcal{L}^\prime_{k,\bm{Z}}(\bm{U_k}, \bm{Z_k})(\bm{\delta Z_k}) &= 0 \quad \forall \bm{\delta  Z_k} \in X_k^{\dG(r)}(V), \\
    \mathcal{L}^\prime_{kh,\bm{Z}}(\bm{U_{\bm{kh}}}, \bm{Z_{\bm{kh}}})(\bm{\delta Z_{\bm{kh}}}) &= 0 \quad \forall \bm{\delta  Z_{\bm{kh}}} \in X_k^{\dG(r)}(V_h^{s_v/s_p}).
\end{align}
We call these equations the primal problems and their solutions $\bm{U}$, $\bm{U_k}$ and $\bm{U}_{\bm{kh}}$ the primal solutions.
Secondly, the stationary points must also satisfy the equations
\begin{align}
    \mathcal{L}^\prime_{\bm{U}}(\bm{U}, \bm{Z})(\bm{\delta U}) &= 0 \quad \forall \bm{\delta  U} \in X(V), \\
    \mathcal{L}^\prime_{k,\bm{U}}(\bm{U_k}, \bm{Z_k})(\bm{\delta U_k}) &= 0 \quad \forall \bm{\delta  U_k} \in X_k^{\dG(r)}(V), \\
    \mathcal{L}^\prime_{kh,\bm{U}}(\bm{U_{\bm{kh}}}, \bm{Z_{\bm{kh}}})(\bm{\delta U_{\bm{kh}}}) &= 0 \quad \forall \bm{\delta  U_{\bm{kh}}} \in X_k^{\dG(r)}(V_h^{s_v/s_p}).
\end{align}
These equations are being called the adjoint or dual problems and their solutions $\bm{Z}$, $\bm{Z_k}$ and $\bm{Z}_{\bm{kh}}$ are the adjoint solutions.
To be able to decide whether the spatial or temporal mesh should be refined, we split the total discretization error into
\begin{align}
    J(\bm{U}) - J(\bm{U_{kh}}) = (J(\bm{U}) - J(\bm{U_k})) + (J(\bm{U_k}) - J(\bm{U_{\bm{kh}}})) \approx \eta_k + \eta_h,
\end{align}
where $\eta_k$ stands for the temporal error estimate and $\eta_h$ stands for the spatial error estimate.
Following Theorem 5.2 in \cite{BeRa12} and neglecting the remainder terms, which are of higher order, we get the error estimates
\begin{equation}\label{eq:general_error_estimator}
\begin{aligned}
    J(\bm{U}) - J(\bm{U_k}) &\approx \frac{1}{2} \left( \rho(\bm{U_k})(\bm{Z}-\bm{Z_k}) + \rho^*(\bm{U_k}, \bm{Z_k})(\bm{U}-\bm{U_k}) \right) =: \eta_k, \\
    J(\bm{U_k}) - J(\bm{U_{\bm{kh}}}) &\approx \frac{1}{2} \left( \rho(\bm{U_{\bm{kh}}})(\bm{Z_k}-\bm{Z_{\bm{kh}}}) + \rho^*(\bm{U_{\bm{kh}}}, \bm{Z_{\bm{kh}}})(\bm{U_k}-\bm{U_{\bm{kh}}}) \right) =: \eta_h,
\end{aligned}
\end{equation}
where we denote the primal and adjoint residuals as
\begin{align}
    \rho(\bm{U})(\bm{\delta Z}) := \mathcal{L}^\prime_{kh,\bm{Z}}(\bm{U}, \bm{Z})(\bm{\delta Z}), \qquad 
    \rho^*(\bm{U},\bm{Z})(\bm{\delta U}) := \mathcal{L}^\prime_{kh,\bm{U}}(\bm{U}, \bm{Z})(\bm{\delta U}) .
\end{align}
We now have two problems to solve: the primal and the adjoint problem. All the remaining variables in the error estimator can be computed via higher order interpolation or via interpolation to a lower order space.

\subsection{Primal problem}
To derive the primal problem, we need to compute the G\^{a}teaux derivatives of the Lagrange functionals $\mathcal{L}$, $\mathcal{L}_{k}$ and $\mathcal{L}_{kh}$ with respect to the adjoint solution $\bm{Z}$. Using the definition of the continuous ($A$), the time-discrete ($A_k$) and the fully discrete ($A_{kh}$) variational formulations, we observe that they are linear in the test functions and hence we have
\begin{align}
    \mathcal{L}^\prime_{\bm{Z}}(\bm{U}, \bm{Z})(\bm{\delta Z}) = -A(\bm{U})(\bm{\delta Z}) &= 0 \quad \forall \bm{\delta  Z} \in X(V), \\
    \mathcal{L}^\prime_{k,\bm{Z}}(\bm{U_k}, \bm{Z_k})(\bm{\delta Z_k}) = -A_k(\bm{U_k})(\bm{\delta Z_k}) &= 0 \quad \forall \bm{\delta  Z_k} \in X_k^{\dG(r)}(V), \\
    \mathcal{L}^\prime_{kh,\bm{Z}}(\bm{U_{\bm{kh}}}, \bm{Z_{\bm{kh}}})(\bm{\delta Z_{\bm{kh}}}) = -A_{kh}(\bm{U_{\bm{kh}}})(\bm{\delta Z_{\bm{kh}}}) &= 0 \quad \forall \bm{\delta  Z_{\bm{kh}}} \in X_k^{\dG(r)}(V_h^{s_v/s_p}).
\end{align}
Therefore, to get the primal solution, we simply need to solve the weak formulation of the Navier-Stokes problem.

\subsection{Adjoint problem}

To derive the adjoint problem, we need to compute the G\^{a}teaux derivatives of the Lagrange functionals $\mathcal{L}$, $\mathcal{L}_{k}$ and $\mathcal{L}_{kh}$ with respect to the primal solution $\bm{U}$. We then get
{\footnotesize
\begin{align}
    \mathcal{L}^\prime_{\bm{U}}(\bm{U}, \bm{Z})(\bm{\delta U}) = J^\prime_{\bm{U}}(\bm{U})(\bm{\delta  U}) -A^\prime_{\bm{U}}(\bm{U})(\bm{\delta U}, \bm{Z}) &= 0 \quad \forall \bm{\delta  U} \in X(V), \\
    \mathcal{L}^\prime_{k,\bm{U}}(\bm{U_k}, \bm{Z_k})(\bm{\delta U_k}) = J^\prime_{\bm{U}}(\bm{U_k})(\bm{\delta  U_k})-A^\prime_{k,\bm{U}}(\bm{U_k})(\bm{\delta U_k}, \bm{Z_k}) &= 0 \quad \forall \bm{\delta  U_k} \in X_k^{\dG(r)}(V), \\
    \mathcal{L}^\prime_{kh,\bm{U}}(\bm{U_{\bm{kh}}}, \bm{Z_{\bm{kh}}})(\bm{\delta U_{\bm{kh}}}) = J^\prime_{\bm{U}}(\bm{U_{\bm{kh}}})(\bm{\delta  U_{\bm{kh}}})-A^\prime_{kh,\bm{U}}(\bm{U_{\bm{kh}}})(\bm{\delta U_{\bm{kh}}}, \bm{Z_{\bm{kh}}}) &= 0 \quad \forall \bm{\delta  U_{\bm{kh}}} \in X_k^{\dG(r)}(V_h^{s_v/s_p}).
\end{align}
}
Therefore, to obtain the adjoint solution, we need to solve an additional equation, the adjoint problem
\begin{align}\label{eq:space_time_discrete_general_adjoint}
    A^\prime_{kh,\bm{U}}(\bm{U_{\bm{kh}}})(\bm{\delta U_{\bm{kh}}}, \bm{Z_{\bm{kh}}}) = J^\prime_{\bm{U}}(\bm{U_{\bm{kh}}})(\bm{\delta  U_{\bm{kh}}}).
\end{align}

\begin{remark}\label{remark:adjoint_linear_problem}
    For linear PDEs, like the Stokes equations, and linear goal functionals, like the mean drag functional, the adjoint problems simplify to
    \begin{align}
        \mathcal{L}^\prime_{\bm{U}}(\bm{U}, \bm{Z})(\bm{\delta U}) = J(\bm{\delta  U}) -A(\bm{\delta U})(\bm{Z}) &= 0 \quad \forall \bm{\delta  U} \in X(V), \\
        \mathcal{L}^\prime_{k,\bm{U}}(\bm{U_k}, \bm{Z_k})(\bm{\delta U_k}) = J(\bm{\delta  U_k})-A_k(\bm{\delta U_k})(\bm{Z_k}) &= 0 \quad \forall \bm{\delta  U_k} \in X_k^{\dG(r)}(V), \\
        \mathcal{L}^\prime_{kh,\bm{U}}(\bm{U_{\bm{kh}}}, \bm{Z_{\bm{kh}}})(\bm{\delta U_{\bm{kh}}}) = J(\bm{\delta  U_{\bm{kh}}})-A_{kh}(\bm{\delta U_{\bm{kh}}})( \bm{Z_{\bm{kh}}}) &= 0 \quad \forall \bm{\delta  U_{\bm{kh}}} \in X_k^{\dG(r)}(V_h^{s_v/s_p}).
    \end{align}
    In this scenario, we have the error identities
    \begin{equation}\label{eq:linear_problem_error_identity}
    \begin{aligned}
        J(\bm{U}) - J(\bm{U_k}) &= -A_k(\bm{U_k})(\bm{Z} - \bm{Z_k}), \\
        J(\bm{U_k}) - J(\bm{U_{\bm{kh}}}) &= -A_{kh}(\bm{U_{\bm{kh}}})(\bm{Z_k} - \bm{Z_{\bm{kh}}}).
    \end{aligned}
    \end{equation}
    We show this error identity for the spatial error contribution.
    Using the linearity of the goal functional and the definition of the adjoint problem, we have
    \begin{align}
        J(\bm{U_k}) - J(\bm{U}_{\bm{kh}}) &= J(\bm{U_k}-\bm{U}_{\bm{kh}}) = A_{k}(\bm{U_k}-\bm{U}_{\bm{kh}})(\bm{Z_k})  = A_{kh}(\bm{U_k}-\bm{U}_{\bm{kh}})(\bm{Z_k}) \\
        &= -A_{kh}(\bm{U}_{\bm{kh}})(\bm{Z_k}) = -A_{kh}(\bm{U}_{\bm{kh}})(\bm{Z_k}-\bm{Z}_{\bm{kh}}).
    \end{align}
    For the temporal error identity we can follow the same steps with a few minor modifications. In the proof of the spatial error, we used that the semilinear form of the time discrete and space-time discrete formulation are the same, i.e. $A_k(\cdot)(\cdot) = A_{kh}(\cdot)(\cdot)$. This does not hold anymore for the continuous and the time discrete formulations, i.e. $A(\cdot)(\cdot) \neq A_k(\cdot)(\cdot)$, since the continuous semilinear form does not include the jump terms. Therefore, we need to replace the continuous semilinear form $A$ by a semilinear form $\tilde{A}$, which contains the jump terms. This replacement is valid, since the continuous solution $\bm{U}$ does not have any jumps and hence $A(\bm{U})(\cdot) = \tilde{A}(\bm{U})(\cdot)$. Finally, we have nonconformity of the continuous and time discrete function spaces, i.e. $X_k^{\dG(r)}(V) \not\subset X(V)$. Hence, for the theory of the error estimator we need to define the continuous variational formulation on the larger space $X(V) + X_k^{\dG(r)}(V)$.
\end{remark}
\noindent By Remark \ref{remark:adjoint_linear_problem}, we can thus write the adjoint problem (\ref{eq:space_time_discrete_general_adjoint}) for the Stokes problem as
\begin{align}
    &\hspace{3cm}A_{kh}(\bm{\delta U_{\bm{kh}}})(\bm{Z_{\bm{kh}}}) = J(\bm{\delta  U_{\bm{kh}}}) \\
    &\Leftrightarrow \sum_{m = 1}^M \int_{I_m}(\partial_t \bm{\delta v_{\bm{kh}}}, \bm{Z_{\bm{kh}}}^v) - (\delta p_{kh}, \nabla_x \cdot \bm{Z_{\bm{kh}}}^v) + \nu (\nabla_x \bm{\delta v_{\bm{kh}}}, \nabla_x \bm{Z_{\bm{kh}}}^v) + (\nabla_x \cdot \bm{\delta v_{\bm{kh}}}, Z_{kh}^p)\ \mathrm{d}t \\
    &\hspace{1cm}+ \sum_{m=0}^{M-1} ([\bm{\delta v_{\bm{kh}}}]_m, \bm{Z}_{\bm{kh},m}^{v,+})  + (\bm{\delta v}_{\bm{kh},0}^-(0) - \bm{\delta v}^0, \bm{Z}_{\bm{kh},0}^{v,-}(0)) = J(\bm{\delta  U_{\bm{kh}}}).
\end{align}
To move the time derivative from the test function to the adjoint solution, we apply integration by parts in time and get
\begin{align}
    &\sum_{m = 1}^M \int_{I_m}(\bm{\delta v_{\bm{kh}}}, -\partial_t \bm{Z_{\bm{kh}}}^v) - (\delta p_{kh}, \nabla_x \cdot \bm{Z_{\bm{kh}}}^v) + \nu (\nabla_x \bm{\delta v_{\bm{kh}}}, \nabla_x \bm{Z_{\bm{kh}}}^v) + (\nabla_x \cdot \bm{\delta v_{\bm{kh}}}, Z_{kh}^p)\ \mathrm{d}t \\
    &\hspace{1cm}- \sum_{m=1}^{M} (\bm{\delta v^-_{\bm{kh},m}}, [\bm{Z}_{\bm{kh}}^{v}]_m)  + (\bm{\delta v}_{\bm{kh},M}^-(T), \bm{Z}_{\bm{kh},M}^{v,-}(T)) = J(\bm{\delta  U_{\bm{kh}}}).
\end{align}
Similarly, to solve the adjoint problem of the Navier-Stokes equations, we now need to solve the general adjoint problem (\ref{eq:space_time_discrete_general_adjoint}), which due to the linearity of our goal functional can be simplified to
\begin{align}
    A^\prime_{kh,\bm{U}}(\bm{U_{\bm{kh}}})(\bm{\delta U_{\bm{kh}}}, \bm{Z_{\bm{kh}}}) = J(\bm{\delta  U_{\bm{kh}}}).
\end{align}
Note that we still need to access the primal solution $\bm{U_{\bm{kh}}}$ to be able to solve the adjoint equation.
The adjoint problem for the Navier-Stokes equations reads
\begin{align}
    &\sum_{m = 1}^M \int_{I_m}(\bm{\delta v_{\bm{kh}}}, -\partial_t \bm{Z_{\bm{kh}}}^v) - (\delta p_{kh}, \nabla_x \cdot \bm{Z_{\bm{kh}}}^v) + \nu (\nabla_x \bm{\delta v_{\bm{kh}}}, \nabla_x \bm{Z_{\bm{kh}}}^v)\ \mathrm{d}t \\
    &\hspace{1cm} +\sum_{m = 1}^M \int_{I_m} ((\bm{\delta v_{\bm{kh}}} \cdot \nabla_x)\bm{v_{\bm{kh}}} + (\bm{v_{\bm{kh}}} \cdot \nabla_x)\bm{\delta v_{\bm{kh}}}, \bm Z^v_{\bm{kh}}) + (\nabla_x \cdot \bm{\delta v_{\bm{kh}}}, Z_{kh}^p)\ \mathrm{d}t\\
    &\hspace{1cm}- \sum_{m=1}^{M} (\bm{\delta v^-_{\bm{kh},m}}, [\bm{Z}_{\bm{kh}}^{v}]_m)  + (\bm{\delta v}_{\bm{kh},M}^-(T), \bm{Z}_{\bm{kh},M}^{v,-}(T)) = J(\bm{\delta  U_{\bm{kh}}}).
\end{align}

\subsection{PU-DWR Navier-Stokes error estimator}
For linear PDEs, like the Stokes problem, and linear goal functionals, we have the error identity 
\begin{align}
    J(\bm{U}) - J(\bm{U_{kh}}) &= [ J(\bm{U}) - J(\bm{U_k}) ] + [ J(\bm{U_k}) - J(\bm{U_{\bm{kh}}}) ] \\
    &= \underbrace{-A_k(\bm{U_k})(\bm{Z} - \bm{Z_k})}_{=: \eta_k} \underbrace{-A_{kh}(\bm{U_{\bm{kh}}})(\bm{Z_k} - \bm{Z_{\bm{kh}}})}_{=: \eta_h} =: \eta,
\end{align}
which was derived in (\ref{eq:linear_problem_error_identity}).
For the nonlinear Navier-Stokes equations, we do not have the error identity (\ref{eq:linear_problem_error_identity}), but need to use approximation (\ref{eq:general_error_estimator}) which contains the primal residual $\rho$ and the adjoint residual $\rho^{\ast}$. However, for mildly nonlinear problems, like incompressible Navier-Stokes flow in the laminar regime, it has been shown that using only the primal residual already produces decent results, see e.g. \cite{becker_rannacher_2001}[Sec. 8] or \cite{Braack2006}. Hence, for the estimation of the temporal and spatial error contributions, we use
\begin{align}
    J(\bm{U}) - J(\bm{U_k}) &\approx -A_k(\bm{U_k})(\bm{Z} - \bm{Z_k}), \\
    J(\bm{U_k}) - J(\bm{U_{\bm{kh}}}) &\approx -A_{kh}(\bm{U_{\bm{kh}}})(\bm{Z_k} - \bm{Z_{\bm{kh}}}).
\end{align}
To localize the spatial and temporal error, we now include the partition of unity $\lbrace \chi_i^m \rbrace_{i,m}$ with $\sum_{m=1}^M \sum_{i \in \mathbb{T}_h^m} \chi_i^m \equiv 1$ in the adjoint weights to localize the error. For the partition of unity we choose linear finite elements in space ($\cG(1)$) and piecewise discontinuous finite elements in time ($\dG(0)$).
This allows for localization of the error to the individual temporal elements and the vertices of the spatial elements. We have 
\begin{align}\label{eq:pu_error_estimators}
    \eta_k = \sum_{m=1}^M \sum_{i \in \mathbb{T}_h^m} \underbrace{-A_k(\bm{U_k})( (\bm{Z}-\bm{Z_k})\chi_i^m)}_{=: \eta_k^{i,m}}, \quad  \eta_h = \sum_{m=1}^M \sum_{i \in \mathbb{T}_h^m} \underbrace{-A_{kh}(\bm{U_{\bm{kh}}})( (\bm{Z_k}-\bm{Z_{\bm{kh}}})\chi_i^m)}_{=: \eta_h^{i,m}}.
\end{align}
Exemplarily, we will show the formula for the computation of the temporal error indicator $\eta_k^{i,m}$ by applying the product rule
\begin{equation}
\small
\begin{aligned}
    \eta_k^{i,m} &= -A_k(\trialc{\bm{U_k}})( \testc{(\bm{Z}-\bm{Z_k})}\puc{\chi_i^m}) \\ 
    &= -\int_{I_m}(\trialc{\partial_t \bm{v_k}}, \testc{(\bm{Z}^v - \bm{Z}_k^v)}\puc{\chi_i^m}) - (\trialc{p_k}, \nabla_x \cdot \lbrace\testc{(\bm{Z}^v - \bm{Z}_k^v)}\puc{\chi_i^m}\rbrace) + \nu (\trialc{\nabla_x \bm{v_k}}, \nabla_x \lbrace\testc{(\bm{Z}^v - \bm{Z}_k^v)}\puc{\chi_i^m}\rbrace) \ \mathrm{d}t \\
    &\quad -\int_{I_m}(\trialc{(\bm{v_{\bm{k}}} \cdot \nabla_x)\bm{ v_{\bm{k}}}}, \testc{(\bm{Z}^v - \bm{Z}_k^v)}\puc{\chi_i^m}) + (\trialc{\nabla_x \cdot \bm{v_k}}, \testc{(Z^p - Z_k^p)}\puc{\chi_i^m})\ \mathrm{d}t - (\trialc{[\bm{v_k}]_m}, \testc{(\bm{Z}_m^{v,+} - \bm{Z}_{k,m}^{v,+})}\puc{\chi_i^m}) \\
    &= -\int_{I_m}(\trialc{\partial_t \bm{v_k}}, \testc{(\bm{Z}^v - \bm{Z}_k^v)}\puc{\chi_i^m}) - (\trialc{p_k}, \testc{\lbrace\nabla_x \cdot (\bm{Z}^v - \bm{Z}_k^v)\rbrace}\puc{\chi_i^m} + \testc{(\bm{Z}^v - \bm{Z}_k^v)} \cdot \puc{(\nabla_x \chi_i^m)})\ \mathrm{d}t \\
    &\quad -\int_{I_m} \nu (\trialc{\nabla_x \bm{v_k}}, \testc{\lbrace\nabla_x(\bm{Z}^v - \bm{Z}_k^v)\rbrace}\puc{\chi_i^m} + \testc{(\bm{Z}^v - \bm{Z}_k^v)} \otimes \puc{(\nabla_x\chi_i^m)} ) \ \mathrm{d}t \\
    &\quad -\int_{I_m}(\trialc{(\bm{v_{\bm{k}}} \cdot \nabla_x)\bm{ v_{\bm{k}}}}, \testc{(\bm{Z}^v - \bm{Z}_k^v)}\puc{\chi_i^m}) + (\trialc{\nabla_x \cdot \bm{v_k}}, \testc{(Z^p - Z_k^p)}\puc{\chi_i^m})\ \mathrm{d}t - (\trialc{[\bm{v_k}]_m}, \testc{(\bm{Z}_m^{v,+} - \bm{Z}_{k,m}^{v,+})}\puc{\chi_i^m}).
\end{aligned}
\end{equation}
\noindent Here $\bm{v} \otimes \bm{w}$ stands for the tensor product of two vectors $\bm{v}$ and $\bm{w}$, i.e. $(\bm{v} \otimes \bm{w})_{i,j} := v_i w_j$. 
With this partition of unity approach the spatial and temporal error can now be localized to each space-time degree of freedom of the $\cG(1)\dG(0)$ space-time discretization and will later indicate which spatial or temporal elements have the largest error and should be refined. However, to test the quality of our error estimator we will initially work with globally refined space-time meshes and monitor how well the estimated error coincides with the true error by computing the effectivity index
\begin{align}
    \Ieff := \left|\frac{\eta}{J(\bm{U})-J(\bm{U}_{\bm{kh}})}\right| = \left|\frac{\eta_k + \eta_h}{J(\bm{U})-J(\bm{U}_{\bm{kh}})}\right|.
\end{align}
We asymptotically expect $\Ieff \underset{k,h \rightarrow 0}{\longrightarrow} 1$ for problems that satisfy a saturation assumption with a constant in the upper bound that converges to 0 for $k,h \rightarrow 0$ \cite{ ThiWi22_arxiv} (based on the ideas from the spatial case \cite{EndtLaWi20}).

\subsection{Adaptive algorithm}
We will use the adaptive algorithm proposed in \cite{ThiWi22_arxiv}, which is a time-slabbing adaptation of the algorithm proposed in \cite{SchmichVexler2008}
to balance the temporal and the spatial discretization error.

\begin{algorithm}[H]
\caption{MARK and REFINE for the time-slabbing approach}\label{algo:equilibration}
\begin{algorithmic}[1]
\Require indicators, equilibration factor $c>0$ and max. number of intervals per slab $N_{\max}$
\State Calculate global temporal estimator 
$\eta_k= \sum\limits_{m=1}^M \sum\limits_{i\in\mathbb{T}_h^m} \eta_{k}^{i,m}$
\State Calculate global spatial estimator 
$\eta_h=  \sum\limits_{m=1}^M \sum\limits_{i\in\mathbb{T}_h^m} \eta_{h}^{i,m}$
\If{$c|\eta_k| \geq |\eta_h|$}
\State Mark temporal elements $I_m$ for refinement based on chosen strategy
\EndIf
\If{$c|\eta_h| \geq |\eta_k|$}
\For{each slab $Q_n$}
\For{each spatial PU-DoF $i\in\mathbb{T}_h(Q_n)$}
\State Calculate $\eta_{h}^{i,Q_n} 
:= \sum\limits_{m=m_n}^{m_{n+1}-1} \eta_{h}^{i,m}$ 
\EndFor
\State Mark and refine elements in $\mathbb{T}_h(Q_n)$ based on chosen strategy using indicators $\eta_{h}^{i,Q_n}$
\EndFor
\EndIf
\For{each slab $Q_n$}
\For{$m=m_n,\dots,m_{n+1}-1$}
\If{$I_m$ is marked}
\State Split/refine $I_m$ into two intervals
\EndIf
\EndFor
\If {$(m_{n+1}-1)-m_n>N_{\max}$}
\State Split slab $Q_n$ into two slabs
\EndIf
\EndFor
\end{algorithmic}
\end{algorithm}
\noindent In our first numerical experiments, we chose a sufficiently big equilibration factor $c$ ($c = 10^7$) to always refine both in space and time. In the future, we want to use a smaller equilibration constant $c$ to only refine either the temporal or the spatial discretization which dominates the total error. 

After having computed the temporal error indicators for each temporal element and the spatial error indicators for each spatial element for each slab, we mark spatial and temporal elements for refinement according to Algorithm \ref{algo:equilibration}. 

\subsection{Evaluation of the error indicators}
Until now, we have not discussed the practical realization of the PU-DWR error estimator (\ref{eq:pu_error_estimators}). To make the error indicators computable, we need to make some approximations and get the unknown solutions by post-processing of the primal and adjoint solutions. In the works of Besier (and co-workers) \cite{SchmichVexler2008, Schmich2009, BeRa12}, the linearization point in the temporal error has been replaced by the primal solution, i.e.
\begin{align}\label{eq:replace_linearization}
    \eta_k = -A_k(\bm{U_k})( \bm{Z} - \bm{Z_k}) \approx -A_k(\bm{U_{\bm{kh}}})( \bm{Z} - \bm{Z_k})
\end{align}
and he has proven that this additional approximation error is small in comparison to the total discretization error \cite{SchmichVexler2008}[Remark 3.2]. Therefore, we will use this simplification in the same way. 
We still need to think about the approximation of the weights $\bm{Z} - \bm{Z_k}$ and $\bm{Z_k} - \bm{Z_{\bm{kh}}}$. This approximation depends on the choice of our finite element space for the adjoint problem and in this work we consider only a mixed order error estimator.

For mixed order, the adjoint solution is of higher order in space and time than the primal solution. More concretely, we choose the primal solution to be $\mathcal{Q}_{2}/\mathcal{Q}_{1}$ in space and $\dG(1)$ in time, whereas the adjoint solution is $\mathcal{Q}_{4}/\mathcal{Q}_{2}$ in space and $\dG(2)$ in time. This setting has already been discussed in \cite{BauBruKoe2021}[Section 4] and in this work we approximate the adjoint weights $\bm{Z} - \bm{Z_k}$ and $\bm{Z_k} - \bm{Z_{\bm{kh}}}$ using the adjoint solution $\bm{Z^{(2)}_{\bm{kh}}} \in X_k^{\dG(2)}\left(V_h^{4/2}\right)$ by 
\begin{align}
    \bm{Z} - \bm{Z_k} &\approx \left(\operatorname{id} -  I_k^{\dG(1)}\right)\bm{Z^{(2)}_{\bm{kh}}}, \\
    \bm{Z_k} - \bm{Z_{\bm{kh}}} &\approx \left(\operatorname{id} - I_h^{2/1}\right)I_k^{\dG(1)}\bm{Z^{(2)}_{\bm{kh}}}.
\end{align}
Here, $I_h^{2/1}: X(V_h^{4/2}) \rightarrow X(V_h^{2/1})$ denotes the interpolation in space from the enriched spatial function space to the primal spatial function space and $I_k^{\dG(1)}: X_k^{\dG(2)}(V) \rightarrow X_k^{\dG(1)}(V)$ denotes the interpolation in time from the enriched temporal function space to the primal temporal function space.

From an implementational point of view this approach is very convenient, since it only requires interpolation of the adjoint solution to a lower order function in space.
This can be done by evaluating the adjoint solution at the space-time degrees of freedom of the lower order space.

%%%%%%%%%%%%%%%%%%%%%%%%%%%%%%%%%%%%%%%%%%%%%%%%%%%%%%%%%%%%
%%                 NUMERICAL TESTS                        %%
%%%%%%%%%%%%%%%%%%%%%%%%%%%%%%%%%%%%%%%%%%%%%%%%%%%%%%%%%%%%
\section{Numerical tests}
\label{sec_tests}

In this section, we perform some numerical experiments to verify the validity of our space-time error estimation framework. In the first two experiments, we consider two-dimensional time-dependent laminar flow around a cylinder, which is being modeled by the Stokes or Navier-Stokes equations. In a further experiment, we model the two-dimensional time-dependent laminar backward-facing step problem 
and do adaptive refinement using a nonlinear functional.
The code is written in the C++ based FEM library deal.II \cite{dealII94, dealii2019design}. 

\subsection{Stokes flow around a cylinder}
\label{sec:stokes}
In the first numerical experiment, we consider the test case 2D-3\footnote{\scriptsize\url{http://www.mathematik.tu-dortmund.de/~featflow/en/benchmarks/cfdbenchmarking/flow/dfg_benchmark3_re100.html}} from the Schäfer/Turek Navier-Stokes benchmarks  \cite{SchaeferTurek1996}.
The 2D-3 benchmark problem describes the two dimensional laminar flow around a circular cylinder with a time-dependent inflow condition. In the following figure, the domain $\Omega$ and the different boundary components $\Gamma_{\text{in}}$, $\Gamma_{\text{out}}$, $\Gamma_{\text{wall}}$ and $\Gamma_{\text{circle}}$ are depicted. \\

\begin{figure}[H]
    \begin{center}
    \includegraphics{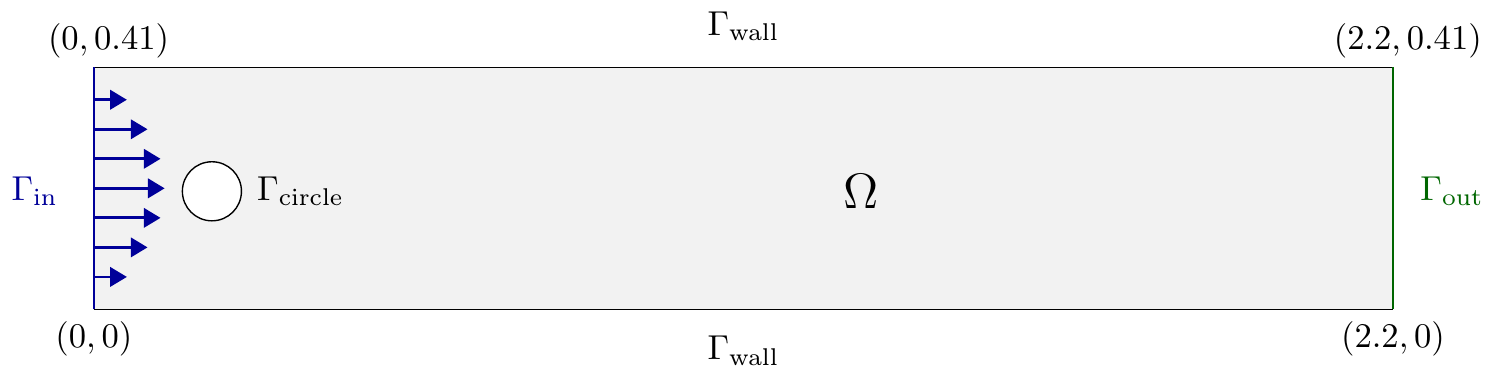}
    \caption{Domain of the 2D Navier-Stokes benchmark problem}
    \label{fig:channel_domain}
    \end{center}
\end{figure}
\noindent Here, the kinematic viscosity is $\nu = 10^{-3}$. The domain is defined as $\Omega := (0,2.2) \times (0,0.41) \setminus B_r(0.2,0.2)$ with $r = 0.05$, where $B_r(x,y)$ denotes the ball around $(x,y)$ with radius $r$ and we denote the time interval by $I := (0, 8)$. On the boundary $\partial \Omega$ we prescribe Dirichlet and Neumann boundary conditions. We apply the no-slip condition, i.e. homogeneous Dirichlet boundary conditions $\bm{v} = \bm{0}$ on $\Gamma_{\text{wall}} \cup \Gamma_{\text{circle}}$, where ${\Gamma_{\text{wall}} = (0,2.2) \times \lbrace 0 \rbrace \cup (0,2.2) \times \lbrace 0.41 \rbrace}$ and ${\Gamma_{\text{circle}} = \partial B_r(0.2,0.2)}$. On the boundary ${\Gamma_{\text{in}} = \lbrace 0 \rbrace \times (0,0.41)}$, we enforce a time-dependent parabolic inflow profile, which is a type of inhomogeneous Dirichlet boundary condition, i.e. $\bm{v} = \bm{v}_D$ on $\Gamma_{\text{in}}$. For the 2D-3 benchmark the inflow parabola is given by
\begin{align}
	\bm{v}_D(t,0,y) = \begin{pmatrix}
		\sin\left(\frac{\pi t}{8}\right)\frac{6y(0.41-y)}{0.41^2} \\
		0
	\end{pmatrix}.
\end{align}
For brevity, we denote the Dirichlet boundary as $\Gamma_D := \Gamma_{\text{in}} \cup \Gamma_{\text{wall}} \cup \Gamma_{\text{circle}}$.
The outflow boundary $\Gamma_{\text{out}} = \lbrace 2.2 \rbrace \times (0,0.41)$ has Neumann boundary conditions. To get a correct outflow profile of the pressure, it has been discussed in \cite{Heywood_Rannacher_Turek_1996} that the (modified) do-nothing condition $\nu (\nabla_x \bm{v}) \cdot \bm{n} - p \bm{n} = \bm{0}$ needs to be enforced on $\Gamma_{\text{out}}$. For the non-symmetric stress tensor this corresponds to homogeneous Neumann boundary conditions, i.e. $\bm\sigma \cdot \bm{n} = \bm{0}$ on $\Gamma_{\text{out}}$.
It can be shown that for this problem the Reynolds number satisfies $0 \leq \operatorname{Re}(t) \leq 100$. 

At first, we consider a simplification of this benchmark problem, where the nonlinearity is being omitted such that we solve the Stokes equations. As mentioned in Remark \ref{remark:adjoint_linear_problem}, this results in an error identity for our primal error estimator and we expect $\Ieff \rightarrow 1$ under  uniform refinement of the space-time meshes. For our calculations, we use the mean drag goal functional
\begin{align}
    J_{\text{drag}}(\bm{U}) := \frac{1}{8} \int_0^8 c_D(t)\ \mathrm{d}t = \frac{20}{8} \int_0^8\int_{\Gamma_{\text{circle}}}\bm\sigma(\bm{U}) \cdot \bm{n} \cdot \bm{e_1}\ \mathrm{d}s\ \mathrm{d}t \approx 0.40284197485629031,
\end{align}
where $\bm{e_1} = (1, 0)^T$, $c_D$ denotes the drag coefficient and the reference value has been calculated on a fine spatio-temporal grid with 10,240 temporal and 1,480,000 spatial degrees of freedom. \\

\subsubsection{Uniform refinement}
\noindent Using our mixed order error estimator, we perform a grid search over the number of primal temporal degrees of freedom $N_k \in \lbrace 20, 40, 80, 160, 320, 640 \rbrace$ and over the number of primal spatial degrees of freedom $N_h \in \lbrace 445, 1610, 6100, 23720, 93520 \rbrace$. Furthermore, we compare Gauss-Legendre and Gauss-Lobatto quadrature for the temporal support points of the $\dG(r)$ method.\\

\noindent\textbf{Gauss-Legendre support points in time}

\begin{table}[H]
    \centering
    \resizebox{0.45\columnwidth}{!}{%
    \begin{tabular}{|r||c|c|c|c|c|c|}
      \hline
     \diagbox{$N_h$}{$N_k$} & 20 & 40 & 80 & 160 & 320 \\ \hline
         445 & 0.61 & 0.61 & 0.61 & 0.61 & 0.61\\ \hline
         1,610 & 0.68 & 0.68 & 0.68 & 0.68 & 0.68\\ \hline
         6,100 & 0.73 & 0.72 & 0.72 & 0.72 & 0.72\\ \hline
         23,720 & 0.77 & 0.75 & 0.75 & 0.75 & 0.75\\ \hline
         93,520 & 0.90 & 0.81 & 0.80 & 0.80 & 0.80\\ \hline
    \end{tabular}
     }
    \caption{Effectivity indices of mixed order for Stokes 2D-3 with \textit{Gauss-Legendre} support points in time}
    \label{table:I_eff_grid_search_stokes}
\end{table}

\noindent We observe in  Table \ref{table:I_eff_grid_search_stokes} that the effectivity index for the Stokes test case is almost constant under temporal refinement. This can be explained by Table \ref{table:eta_k_grid_search_stokes} and Table \ref{table:eta_h_grid_search_stokes}, where we find that the discretization error consists mainly of the spatial discretization error while the temporal discretization error is a few orders of magnitude smaller for this test case. Moreover, it can be seen that the effectivity index converges towards $1$ for uniform spatial refinement, which proves numerically that in the limit our PU-DWR error estimator and the true error coincide.

\begin{table}[H]
    \centering
    \resizebox{0.7\columnwidth}{!}{%
    \begin{tabular}{|r||c|c|c|c|c|c|}
      \hline
     \diagbox{$N_h$}{$N_k$} & 20 & 40 & 80 & 160 & 320 \\ \hline
         445 & \cellcolor[rgb]{1,1,1}$-2.7901 \cdot 10^{-6}$ & \cellcolor[rgb]{0.9,0.9,0.9}$-5.2875 \cdot 10^{-7}$ & \cellcolor[rgb]{1,1,1}$-1.0091 \cdot 10^{-7}$ & \cellcolor[rgb]{0.9,0.9,0.9}$-1.7218 \cdot 10^{-8}$ & \cellcolor[rgb]{1,1,1}$-2.5813 \cdot 10^{-9}$\\ \hline
         1,610 & \cellcolor[rgb]{1,1,1}$-1.0153 \cdot 10^{-5}$ & \cellcolor[rgb]{0.9,0.9,0.9}$-1.5887 \cdot 10^{-6}$ & \cellcolor[rgb]{1,1,1}$-2.2938 \cdot 10^{-7}$ & \cellcolor[rgb]{0.9,0.9,0.9}$-3.1912 \cdot 10^{-8}$ & \cellcolor[rgb]{1,1,1}$-4.3840 \cdot 10^{-9}$\\ \hline
         6,100 & \cellcolor[rgb]{1,1,1}$-1.1674 \cdot 10^{-5}$ & \cellcolor[rgb]{0.9,0.9,0.9}$-1.9094 \cdot 10^{-6}$ & \cellcolor[rgb]{1,1,1}$-2.8545 \cdot 10^{-7}$ & \cellcolor[rgb]{0.9,0.9,0.9}$-4.0021 \cdot 10^{-8}$ & \cellcolor[rgb]{1,1,1}$-5.3819 \cdot 10^{-9}$\\ \hline
         23,720 & \cellcolor[rgb]{1,1,1}$-1.1797 \cdot 10^{-5}$ & \cellcolor[rgb]{0.9,0.9,0.9}$-1.9447 \cdot 10^{-6}$ & \cellcolor[rgb]{1,1,1}$-2.9447 \cdot 10^{-7}$ & \cellcolor[rgb]{0.9,0.9,0.9}$-4.1974 \cdot 10^{-8}$ & \cellcolor[rgb]{1,1,1}$-5.7309 \cdot 10^{-9}$\\ \hline
         93,520 & \cellcolor[rgb]{1,1,1}$-1.1806 \cdot 10^{-5}$ & \cellcolor[rgb]{0.9,0.9,0.9}$-1.9471 \cdot 10^{-6}$ & \cellcolor[rgb]{1,1,1}$-2.9515 \cdot 10^{-7}$ & \cellcolor[rgb]{0.9,0.9,0.9}$-4.2178 \cdot 10^{-8}$ & \cellcolor[rgb]{1,1,1}$-5.7855 \cdot 10^{-9}$\\ \hline
    \end{tabular}
     }
    \caption{Temporal error of mixed order for Stokes 2D-3 with \textit{Gauss-Legendre} support points in time}
    \label{table:eta_k_grid_search_stokes}
\end{table}

\noindent Looking at the columns of Table \ref{table:eta_k_grid_search_stokes}, i.e. keeping the number of temporal degrees of freedom fixed and only refining uniformly in space, we see that the temporal error remains almost constant. This verifies that our temporal error does not depend on spatial refinement.

\begin{table}[H]
    \centering
    \resizebox{0.7\columnwidth}{!}{%
    \begin{tabular}{|r||c|c|c|c|c|c|}
      \hline
     \diagbox{$N_h$}{$N_k$} & 20 & 40 & 80 & 160 & 320 \\ \hline
         445 & \cellcolor[rgb]{1,1,1}$3.7394 \cdot 10^{-2}$ & \cellcolor[rgb]{1,1,1}$3.7397 \cdot 10^{-2}$ & \cellcolor[rgb]{1,1,1}$3.7398 \cdot 10^{-2}$ & \cellcolor[rgb]{1,1,1}$3.7398 \cdot 10^{-2}$ & \cellcolor[rgb]{1,1,1}$3.7398 \cdot 10^{-2}$\\ \hline
         1,610 & \cellcolor[rgb]{0.9,0.9,0.9}$1.3111 \cdot 10^{-2}$ & \cellcolor[rgb]{0.9,0.9,0.9}$1.3110 \cdot 10^{-2}$ & \cellcolor[rgb]{0.9,0.9,0.9}$1.3110 \cdot 10^{-2}$ & \cellcolor[rgb]{0.9,0.9,0.9}$1.3110 \cdot 10^{-2}$ & \cellcolor[rgb]{0.9,0.9,0.9}$1.3110 \cdot 10^{-2}$\\ \hline
         6,100 & \cellcolor[rgb]{1,1,1}$4.0683 \cdot 10^{-3}$ & \cellcolor[rgb]{1,1,1}$4.0674 \cdot 10^{-3}$ & \cellcolor[rgb]{1,1,1}$4.0673 \cdot 10^{-3}$ & \cellcolor[rgb]{1,1,1}$4.0673 \cdot 10^{-3}$ & \cellcolor[rgb]{1,1,1}$4.0673 \cdot 10^{-3}$\\ \hline
         23,720 & \cellcolor[rgb]{0.9,0.9,0.9}$1.1460 \cdot 10^{-3}$ & \cellcolor[rgb]{0.9,0.9,0.9}$1.1457 \cdot 10^{-3}$ & \cellcolor[rgb]{0.9,0.9,0.9}$1.1457 \cdot 10^{-3}$ & \cellcolor[rgb]{0.9,0.9,0.9}$1.1457 \cdot 10^{-3}$ & \cellcolor[rgb]{0.9,0.9,0.9}$1.1457 \cdot 10^{-3}$\\ \hline
         93,520 & \cellcolor[rgb]{1,1,1}$3.0495 \cdot 10^{-4}$ & \cellcolor[rgb]{1,1,1}$3.0486 \cdot 10^{-4}$ & \cellcolor[rgb]{1,1,1}$3.0485 \cdot 10^{-4}$ & \cellcolor[rgb]{1,1,1}$3.0485 \cdot 10^{-4}$ & \cellcolor[rgb]{1,1,1}$3.0484 \cdot 10^{-4}$\\ \hline
    \end{tabular}
     }
    \caption{Spatial error of mixed order for Stokes 2D-3 with \textit{Gauss-Legendre} support points in time}
    \label{table:eta_h_grid_search_stokes}
\end{table}

\noindent Looking at the rows of Figure \ref{table:eta_h_grid_search_stokes}, i.e. keeping the number of spatial degrees of freedom fixed and only refining uniformly in time, we see that the spatial error remains almost constant. This verifies that our spatial error does not depend on temporal refinement. Additionally, we observe quadratic convergence of the spatial error with respect to the spatial mesh size $h$.\\

\noindent\textbf{Gauss-Lobatto support points in time}

\begin{table}[H]
    \centering
    \resizebox{0.45\columnwidth}{!}{%
    \begin{tabular}{|r||c|c|c|c|c|c|}
      \hline
     \diagbox{$N_h$}{$N_k$} & 20 & 40 & 80 & 160 & 320 \\ \hline
         445 & 0.58 & 0.61 & 0.61 & 0.61 & 0.61\\ \hline
         1,610 & 0.58 & 0.65 & 0.67 & 0.68 & 0.68\\ \hline
         6,100 & 0.45 & 0.63 & 0.70 & 0.72 & 0.72\\ \hline
         23,720 & 0.23 & 0.49 & 0.66 & 0.73 & 0.75\\ \hline
         93,520 & 0.08 & 0.25 & 0.52 & 0.70 & 0.77\\ \hline
    \end{tabular}
     }
    \caption{Effectivity indices of mixed order for Stokes 2D-3 with \textit{Gauss-Lobatto} support points in time}
    \label{table:I_eff_grid_search_stokes_gauss_lobatto}
\end{table}

\noindent In contrast to Table \ref{table:I_eff_grid_search_stokes}, we observe in Table \ref{table:I_eff_grid_search_stokes_gauss_lobatto} that the effectivity indices with Gauss-Lobatto support points in time are smaller for coarse temporal meshes. However, for fine temporal scales the effectivity indices from Table \ref{table:I_eff_grid_search_stokes_gauss_lobatto} converge to the effectivity indices of Table \ref{table:I_eff_grid_search_stokes}.

\begin{table}[H]
    \centering
    \resizebox{0.7\columnwidth}{!}{%
    \begin{tabular}{|r||c|c|c|c|c|c|}
      \hline
     \diagbox{$N_h$}{$N_k$} & 20 & 40 & 80 & 160 & 320 \\ \hline
         445 & \cellcolor[rgb]{1,1,1}$-1.0450 \cdot 10^{-3}$ & \cellcolor[rgb]{0.9,0.9,0.9}$-3.3808 \cdot 10^{-4}$ & \cellcolor[rgb]{1,1,1}$8.4236 \cdot 10^{-5}$ & \cellcolor[rgb]{0.9,0.9,0.9}$2.6087 \cdot 10^{-4}$ & \cellcolor[rgb]{1,1,1}$3.2106 \cdot 10^{-4}$\\ \hline
         1,610 & \cellcolor[rgb]{1,1,1}$-2.4703 \cdot 10^{-4}$ & \cellcolor[rgb]{0.9,0.9,0.9}$-1.6773 \cdot 10^{-4}$ & \cellcolor[rgb]{1,1,1}$-8.1390 \cdot 10^{-5}$ & \cellcolor[rgb]{0.9,0.9,0.9}$-1.2743 \cdot 10^{-5}$ & \cellcolor[rgb]{1,1,1}$2.5248 \cdot 10^{-5}$\\ \hline
         6,100 & \cellcolor[rgb]{1,1,1}$-4.2596 \cdot 10^{-5}$ & \cellcolor[rgb]{0.9,0.9,0.9}$-2.9024 \cdot 10^{-5}$ & \cellcolor[rgb]{1,1,1}$-2.2523 \cdot 10^{-5}$ & \cellcolor[rgb]{0.9,0.9,0.9}$-1.5217 \cdot 10^{-5}$ & \cellcolor[rgb]{1,1,1}$-7.6783 \cdot 10^{-6}$\\ \hline
         23,720 & \cellcolor[rgb]{1,1,1}$-1.7540 \cdot 10^{-5}$ & \cellcolor[rgb]{0.9,0.9,0.9}$-4.5460 \cdot 10^{-6}$ & \cellcolor[rgb]{1,1,1}$-2.5516 \cdot 10^{-6}$ & \cellcolor[rgb]{0.9,0.9,0.9}$-2.0765 \cdot 10^{-6}$ & \cellcolor[rgb]{1,1,1}$-1.6699 \cdot 10^{-6}$\\ \hline
         93,520 & \cellcolor[rgb]{1,1,1}$-1.5636 \cdot 10^{-5}$ & \cellcolor[rgb]{0.9,0.9,0.9}$-2.4541 \cdot 10^{-6}$ & \cellcolor[rgb]{1,1,1}$-4.8653 \cdot 10^{-7}$ & \cellcolor[rgb]{0.9,0.9,0.9}$-2.0288 \cdot 10^{-7}$ & \cellcolor[rgb]{1,1,1}$-1.5701 \cdot 10^{-7}$\\ \hline
    \end{tabular}
     }
    \caption{Temporal error of mixed order for Stokes 2D-3 with \textit{Gauss-Lobatto} support points in time}
    \label{table:eta_k_grid_search_stokes_gauss_lobatto}
\end{table}

\noindent In Table \ref{table:eta_k_grid_search_stokes_gauss_lobatto}, we observe that the temporal error with Gauss-Lobatto support points is larger than in Table \ref{table:eta_k_grid_search_stokes} for Gauss-Legendre. Therefore, from now on we only use Gauss-Legendre support points in time.

\begin{table}[H]
    \centering
    \resizebox{0.7\columnwidth}{!}{%
    \begin{tabular}{|r||c|c|c|c|c|c|}
      \hline
     \diagbox{$N_h$}{$N_k$} & 20 & 40 & 80 & 160 & 320 \\ \hline
         445 & \cellcolor[rgb]{1,1,1}$3.8122 \cdot 10^{-2}$ & \cellcolor[rgb]{1,1,1}$3.7657 \cdot 10^{-2}$ & \cellcolor[rgb]{1,1,1}$3.7294 \cdot 10^{-2}$ & \cellcolor[rgb]{1,1,1}$3.7132 \cdot 10^{-2}$ & \cellcolor[rgb]{1,1,1}$3.7076 \cdot 10^{-2}$\\ \hline
         1,610 & \cellcolor[rgb]{0.9,0.9,0.9}$1.3234 \cdot 10^{-2}$ & \cellcolor[rgb]{0.9,0.9,0.9}$1.3249 \cdot 10^{-2}$ & \cellcolor[rgb]{0.9,0.9,0.9}$1.3185 \cdot 10^{-2}$ & \cellcolor[rgb]{0.9,0.9,0.9}$1.3121 \cdot 10^{-2}$ & \cellcolor[rgb]{0.9,0.9,0.9}$1.3085 \cdot 10^{-2}$\\ \hline
         6,100 & \cellcolor[rgb]{1,1,1}$4.0609 \cdot 10^{-3}$ & \cellcolor[rgb]{1,1,1}$4.0857 \cdot 10^{-3}$ & \cellcolor[rgb]{1,1,1}$4.0874 \cdot 10^{-3}$ & \cellcolor[rgb]{1,1,1}$4.0820 \cdot 10^{-3}$ & \cellcolor[rgb]{1,1,1}$4.0749 \cdot 10^{-3}$\\ \hline
         23,720 & \cellcolor[rgb]{0.9,0.9,0.9}$1.1383 \cdot 10^{-3}$ & \cellcolor[rgb]{0.9,0.9,0.9}$1.1456 \cdot 10^{-3}$ & \cellcolor[rgb]{0.9,0.9,0.9}$1.1473 \cdot 10^{-3}$ & \cellcolor[rgb]{0.9,0.9,0.9}$1.1476 \cdot 10^{-3}$ & \cellcolor[rgb]{0.9,0.9,0.9}$1.1473 \cdot 10^{-3}$\\ \hline
         93,520 & \cellcolor[rgb]{1,1,1}$3.0247 \cdot 10^{-4}$ & \cellcolor[rgb]{1,1,1}$3.0437 \cdot 10^{-4}$ & \cellcolor[rgb]{1,1,1}$3.0485 \cdot 10^{-4}$ & \cellcolor[rgb]{1,1,1}$3.0496 \cdot 10^{-4}$ & \cellcolor[rgb]{1,1,1}$3.0499 \cdot 10^{-4}$\\ \hline
    \end{tabular}
     }
    \caption{Spatial error of mixed order for Stokes 2D-3 with \textit{Gauss-Lobatto} support points in time}
    \label{table:eta_h_grid_search_stokes_gauss_lobatto}
\end{table}

\noindent Finally, in Table \ref{table:eta_h_grid_search_stokes_gauss_lobatto} we get similar spatial error estimates as in  Table \ref{table:eta_h_grid_search_stokes}. The only difference being that the spatial error estimator changes slightly less under temporal refinement for Gauss-Legendre support points in time, which further motivates our preferred usage of this quadrature formula.

\subsubsection{Adaptive refinement}
\noindent We now apply the space-time adaptivity Algorithm \ref{algo:equilibration} starting from the coarsest spatial mesh with $N_h = 445$ primal DoFs and $M = 20$ slabs with one temporal element each, i.e. we have $N_k = 40$ initial temporal degrees of freedom due to the $\dG(1)$ discretization. This results in $17,800$ space-time degrees of freedom.

As a marking strategy in space, we use the averaging strategy, which marks all elements $K \in \mathbb{T}_h$ such that
\begin{align}
    \eta_h\restrict{K} > \frac{\alpha}{|\mathbb{T}_h|}\sum_{K^\prime \in \mathbb{T}_h}\eta_h\restrict{K^\prime}.
\end{align}
In our computations, we used $\alpha = 1.1$.
As a marking strategy in time, we use a fixed rate strategy (fixed number), where the $\varepsilon \%$ of temporal cells with largest error are being refined. In our computations, we used $\varepsilon = 75$.

\begin{table}[H]
    \centering
    \resizebox{\columnwidth}{!}{%
    \begin{tabular}{|r|r|r|c|c|c|r|r|}
      \hline
     \#DoF(primal) & \#DoF(adjoint) & $M$ & $\eta_k$ & $\eta_h$ & $\eta$ & $J(\bm{U}) - J(\bm{U}_{\bm{kh}})$ & $\Ieff$ \\ \hline
          17,800 &  96,600 & 20 & $-5.2875 \cdot 10^{-7}$ & $3.7397 \cdot 10^{-2}$ & $3.7397 \cdot 10^{-2}$ & $6.0872 \cdot 10^{-2}$ & 0.61 \\ \hline
          55,328 &  303,288 & 36 & $5.6502 \cdot 10^{-8}$ & $1.5619 \cdot 10^{-2}$ & $1.5619 \cdot 10^{-2}$ & $2.1453 \cdot 10^{-2}$ & 0.73 \\ \hline
          259,712 &  1,468,032 & 64 & $-3.8869 \cdot 10^{-8}$ & $4.1840 \cdot 10^{-3}$ & $4.1840 \cdot 10^{-3}$ & $5.7545 \cdot 10^{-3}$ & 0.73 \\ \hline
          902,838 &  5,137,110 & 113 & $-8.8754 \cdot 10^{-9}$ & $1.1668 \cdot 10^{-3}$ & $1.1668 \cdot 10^{-3}$ & $1.5598 \cdot 10^{-3}$ & 0.75 \\ \hline
          3,438,496 &  19,668,018 & 199 & $-7.5194 \cdot 10^{-9}$ & $3.1632 \cdot 10^{-4}$ & $3.1632 \cdot 10^{-4}$ & $4.2411 \cdot 10^{-4}$ & 0.75 \\ \hline
    \end{tabular}
    }
    \caption{Adaptive refinement of mixed order on dynamic meshes for Stokes 2D-3}
    \label{table:adaptive_refinement_Stokes}
\end{table}

\begin{figure}[H]%
    \centering
    \subfloat[Drag coefficient]{
        \includegraphics[clip,trim=3 0 4 0, width=0.45\textwidth]{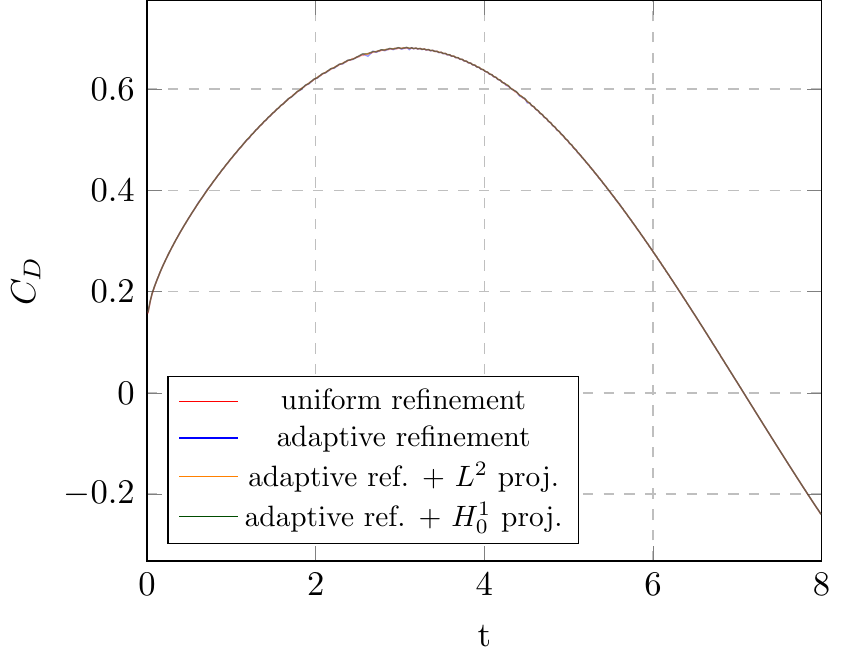}
    }\qquad 
    \subfloat[Lift coefficient]{
        \includegraphics[clip,trim=3 0 4 0,     width=0.45\textwidth]{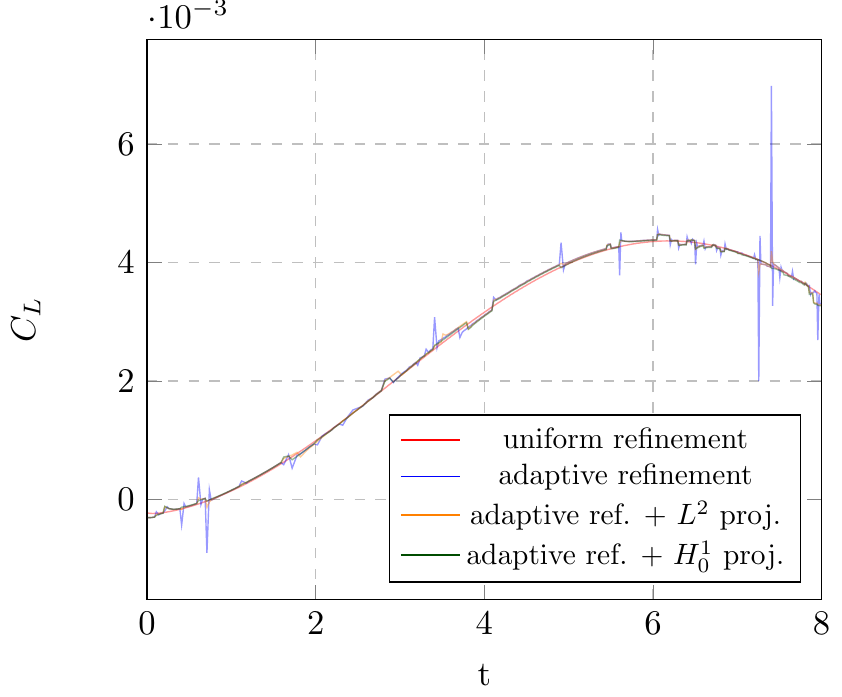}
    }
    \caption{Trajectory of the drag/lift coefficient for the Stokes problem for adaptive refinement}%
    \label{fig:stokes_goal_functionals_adaptive}
\end{figure}

\noindent In Table \ref{table:adaptive_refinement_Stokes}, we show the results from five adaptive refinement loops with a mixed order error estimator for the Stokes 2D-3 problem. We do not use a divergence-free projection here, which leads to oscillations in the goal functionals (cf. Figure \ref{fig:stokes_goal_functionals_adaptive}). However, the error estimates for the Stokes problem are still acceptable and the effectivity index does not change much.

\begin{table}[H]
    \centering
    \resizebox{\columnwidth}{!}{%
    \begin{tabular}{|r|r|r|c|c|c|r|r|}
      \hline
     \#DoF(primal) & \#DoF(adjoint) & $M$ & $\eta_k$ & $\eta_h$ & $\eta$ & $J(\bm{U}) - J(\bm{U}_{\bm{kh}})$ & $\Ieff$ \\ \hline
          17,800 &  96,600 & 20 & $-1.7268 \cdot 10^{-6}$ & $3.7635 \cdot 10^{-2}$ & $3.7634 \cdot 10^{-2}$ & $5.9702 \cdot 10^{-2}$ & 0.63 \\ \hline
          55,092 &  301,932 & 36 & $-4.7644 \cdot 10^{-8}$ & $1.5736 \cdot 10^{-2}$ & $1.5736 \cdot 10^{-2}$ & $2.1001 \cdot 10^{-2}$ & 0.75 \\ \hline
          259,712 &  1,468,032 & 64 & $-4.2556 \cdot 10^{-8}$ & $4.2011 \cdot 10^{-3}$ & $4.2011 \cdot 10^{-3}$ & $5.3248 \cdot 10^{-3}$ & 0.79 \\ \hline
          902,838 &  5,137,110 & 113 & $-1.0787 \cdot 10^{-8}$ & $1.1719 \cdot 10^{-3}$ & $1.1718 \cdot 10^{-3}$ & $1.3084 \cdot 10^{-3}$ & 0.90 \\ \hline
          3,437,792 &  19,663,926 & 199 & $-7.5608 \cdot 10^{-9}$ & $3.1791 \cdot 10^{-4}$ & $3.1790 \cdot 10^{-4}$ & $2.7098 \cdot 10^{-4}$ & 1.17 \\ \hline
    \end{tabular}
    }
    \caption{Adaptive refinement of mixed order on dynamic meshes for Stokes 2D-3 with divergence-free $L^2$ projection}
    \label{table:adaptive_refinement_Stokes_L2}
\end{table}
\noindent In Table \ref{table:adaptive_refinement_Stokes_L2}, we have the same setup as in Table \ref{table:adaptive_refinement_Stokes} and additionally employ a divergence-free $L^2$ projection from Section \ref{sec:div_free_projection}. We observe that ensuring the incompressibility of the snapshot on the new mesh yields effectivity indices closer to 1 and the true error in the goal functional between the reference solution and the finite element solution is smaller than in  Table \ref{table:adaptive_refinement_Stokes}.

\begin{table}[H]
    \centering
    \resizebox{\columnwidth}{!}{%
    \begin{tabular}{|r|r|r|c|c|c|r|r|}
      \hline
     \#DoF(primal) & \#DoF(adjoint) & $M$ & $\eta_k$ & $\eta_h$ & $\eta$ & $J(\bm{U}) - J(\bm{U}_{\bm{kh}})$ & $\Ieff$ \\ \hline
          17,800 &  96,600 & 20 & $-1.9364 \cdot 10^{-6}$ & $3.7704 \cdot 10^{-2}$ & $3.7702 \cdot 10^{-2}$ & $5.9309 \cdot 10^{-2}$ & 0.64 \\ \hline
          55,092 &  301,932 & 36 & $-9.5524 \cdot 10^{-8}$ & $1.5745 \cdot 10^{-2}$ & $1.5745 \cdot 10^{-2}$ & $2.0791 \cdot 10^{-2}$ & 0.76 \\ \hline
          259,712 &  1,468,032 & 64 & $-4.7248 \cdot 10^{-8}$ & $4.1966 \cdot 10^{-3}$ & $4.1966 \cdot 10^{-3}$ & $5.0362 \cdot 10^{-3}$ & 0.83 \\ \hline
          902,838 &  5,137,110 & 113 & $-1.1303 \cdot 10^{-8}$ & $1.1696 \cdot 10^{-3}$ & $1.1696 \cdot 10^{-3}$ & $1.1248 \cdot 10^{-3}$ & 1.04 \\ \hline
          3,439,256 &  19,672,314 & 199 & $-7.3748 \cdot 10^{-9}$ & $3.1672 \cdot 10^{-4}$ & $3.1672 \cdot 10^{-4}$ & $1.5370 \cdot 10^{-4}$ & 2.06 \\ \hline
    \end{tabular}
    }
    \caption{Adaptive refinement of mixed order on dynamic meshes for Stokes 2D-3 with divergence-free $H^1_0$ projection}
    \label{table:adaptive_refinement_Stokes_H1}
\end{table}
\noindent In contrast to Table \ref{table:adaptive_refinement_Stokes_L2}, we now use a divergence-free $H^1_0$ projection in Table \ref{table:adaptive_refinement_Stokes_H1}, but we make similar observations when it comes to the effectivity index and the true error.

\begin{figure}[H]%
    \centering
    \includegraphics{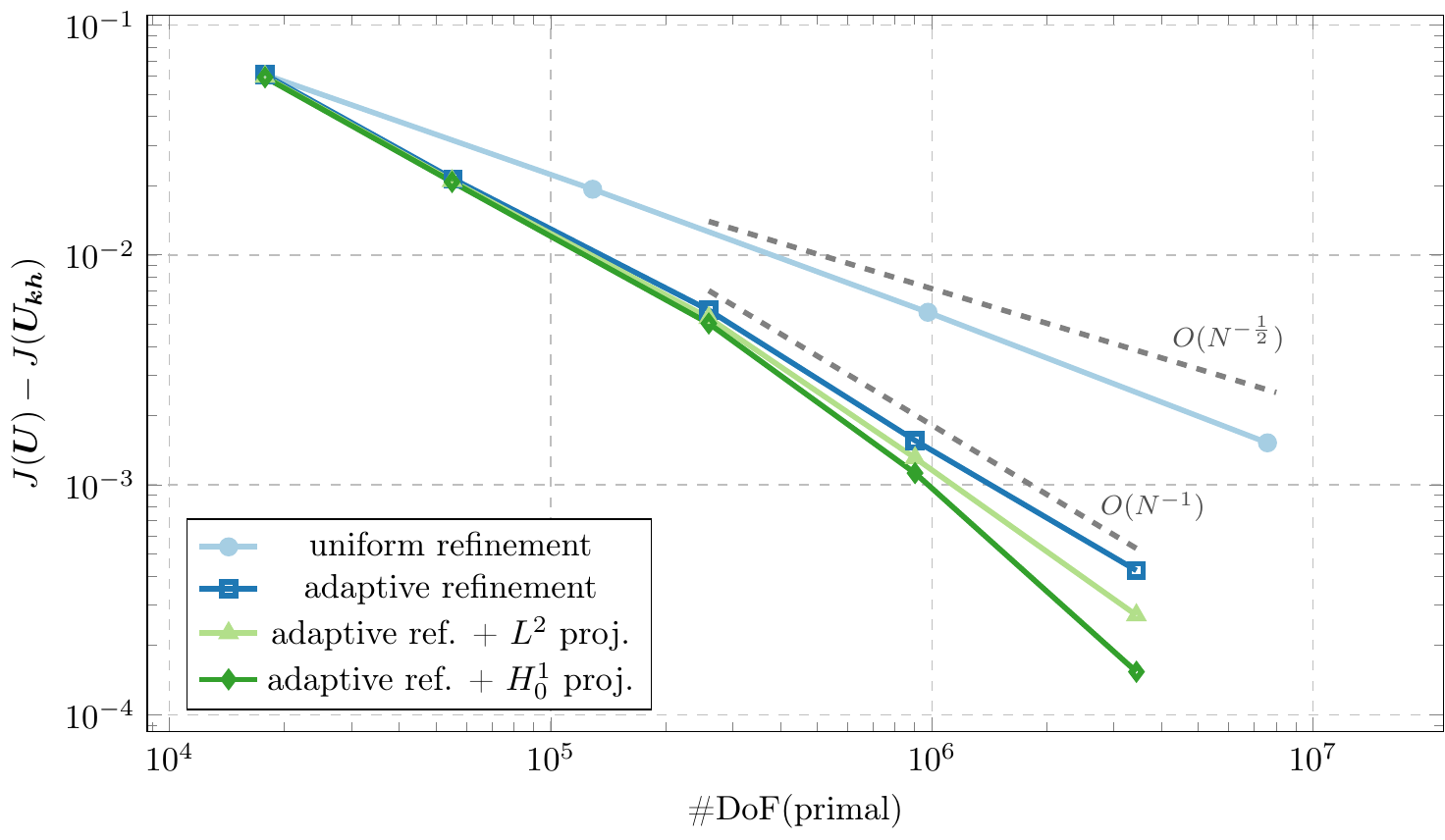}
    \caption{Comparison of the error in the mean drag for the Stokes problem for uniform and adaptive spatio-temporal refinement}%
    \label{fig:stokes_comparison_error_uniform_and_adaptive_refinement}
\end{figure}
\noindent In Figure \ref{fig:stokes_comparison_error_uniform_and_adaptive_refinement}, we start with a common coarse spatio-temporal mesh and iteratively refine uniformly in space and time or adaptively refine in space and time based on Algorithm \ref{algo:equilibration}. As expected the dual weighted residual method driven mesh refinement leads to a faster convergence in the error in the mean drag functional. In particular, for adaptive refinement a lot fewer degrees of freedom are needed to reach a prescribed accuracy in the goal functional, e.g. to achieve a tolerance of $2 \cdot 10^{-3}$ using uniform refinement requires 7,590,400 space-time degrees of freedom, whereas adaptive refinement only needs 902,838 space-time degrees of freedom for the primal problem.

\noindent In Figure \ref{fig:stokes_goal_functionals_adaptive}, we compare the trajectory for the drag and lift coefficients over time for the Stokes 2D-3 problem and the different adaptive refinement methods. 
The drag trajectory is mostly smooth with only minor sawtooth-like oscillations near the maximal drag. 
However, the lift coefficient displays large spikes when using adaptive refinement without ensuring the incompressibility of the solution on the next spatial mesh. 
This underlines the necessity of methods like the divergence-free $L^2$ and $H^1_0$ projections from Besier and Wollner \cite{BesierWollner2012}, which lead to smoother lift coefficient trajectories that more closely match the uniformly refined solution.

\subsection{Navier-Stokes flow around a cylinder}
\label{sec:nse_2d_3}
\noindent We now include the nonlinear term in the variational formulation and solve the Navier-Stokes equations. In contrast to the Stokes problem from Section \ref{sec:stokes}, there are plenty of reference results for the Navier-Stokes 2D-3 benchmark problem in the literature \cite{SchaeferTurek1996, John2004} and we decided to use the publicly available results of the FeatFlow software\footnote{\url{http://www.mathematik.tu-dortmund.de/~featflow/en/benchmarks/cfdbenchmarking/flow/dfg_benchmark3_re100.html}} of the group at the TU Dortmund around Stefan Turek. The FeatFlow authors discretize in space with $Q_2/P_1^{\text{disc}}$ elements and in time they discretize with the Crank Nicholson scheme. 
For their finest calculations, they use a fixed time step size of $k = 0.000625$ and $667,264$ degrees of freedom in space.
We use their mean drag
\begin{align}
    J_{\text{drag}}(\bm{U}) \approx 1.6031368118815639
\end{align}
as our reference solution.
As marking strategies we again used the averaging strategy in space and the 
fixer rate strategy in time with the same parameters as for the 
linear case ($\alpha = 1.1$ and $\varepsilon = 75$).

\begin{table}[H]
    \centering
    \resizebox{\columnwidth}{!}{%
    \begin{tabular}{|r|r|r|c|c|c|r|r|}
      \hline
     \#DoF(primal) & \#DoF(adjoint) & $M$ & $\eta_k$ & $\eta_h$ & $\eta$ & $J(\bm{U}) - J(\bm{U}_{\bm{kh}})$ & $\Ieff$ \\ \hline
          17,800 &  96,600 & 20 & $-8.4213 \cdot 10^{-6}$ & $3.8557 \cdot 10^{-1}$ & $3.8557 \cdot 10^{-1}$ & $5.5869 \cdot 10^{-1}$ & 0.69 \\ \hline
          128,800 &  732,000 & 40 & $1.3424 \cdot 10^{-3}$ & $2.8297 \cdot 10^{-1}$ & $2.8431 \cdot 10^{-1}$ & $2.9690 \cdot 10^{-1}$ & 0.96 \\ \hline
          976,000 &  5,692,800 & 80 & $6.8970 \cdot 10^{-2}$ & $1.3391 \cdot10 ^{-1}$ & $2.0288 \cdot 10^{-1}$ & $1.3978 \cdot 10^{-1}$ & 1.45 \\ \hline
          7,590,400 &  44,889,600 & 160 & $5.5834 \cdot 10^{-3}$ & $2.4517 \cdot 10^{-2}$ & $3.0100 \cdot 10^{-2}$ & $2.5428 \cdot 10^{-2}$ & 1.18 \\ \hline
    \end{tabular}
    }
    \caption{Uniform refinement of mixed order on dynamic meshes for Navier-Stokes 2D-3}
    \label{table:uniform_refinement_nse}
\end{table}

\noindent In Table \ref{table:uniform_refinement_nse}, we observe that the effectivity indices are close to 1 under uniform spatio-temporal refinement for the Navier-Stokes 2D-3 problem. We want to remark again that for this nonlinear problem our primal error estimator does no longer satisfy the error identity, but we can computationally verify its accuracy.

\begin{table}[H]
    \centering
    \resizebox{\columnwidth}{!}{%
    \begin{tabular}{|r|r|r|c|c|c|r|r|}
      \hline
     \#DoF(primal) & \#DoF(adjoint) & $M$ & $\eta_k$ & $\eta_h$ & $\eta$ & $J(\bm{U}) - J(\bm{U}_{\bm{kh}})$ & $\Ieff$ \\ \hline
          17,800 &  96,600 & 20 & $-8.4213 \cdot 10^{-6}$ & $3.8557 \cdot 10^{-1}$ & $3.8557 \cdot 10^{-1}$ & $5.5869 \cdot 10^{-1}$ & 0.69 \\ \hline
          63,336 &  350,244 & 36 & $2.3497 \cdot 10^{-7}$ & $2.6516 \cdot 10^{-1}$ & $2.6516 \cdot 10^{-1}$ & $2.6719 \cdot 10^{-1}$ & 0.99 \\ \hline
          228,568 &  1,285,878 & 64 & $1.7174 \cdot 10^{-3}$ & $-1.0803 \cdot 10^{-2}$ & $9.0853 \cdot 10^{-3}$ & $1.2438 \cdot 10^{-1}$ & 0.07 \\ \hline
          835,942 &  4,749,933 & 113 & $-1.3694 \cdot 10^{-1}$ & $1.7656 \cdot 10^{-1}$ & $3.9619 \cdot 10^{-2}$ & $2.7110 \cdot 10^{-2}$ & 1.46 \\ \hline
          3,008,930 &  17,274,966 & 199 & $3.9451 \cdot 10^{-3}$ & $1.1779 \cdot 10^{-2}$ & $1.5724 \cdot 10^{-2}$ & $1.1074 \cdot 10^{-2}$ & 1.42 \\ \hline
    \end{tabular}
    }
    \caption{Adaptive refinement of mixed order on dynamic meshes for Navier-Stokes 2D-3}
    \label{table:adaptive_refinement_nse}
\end{table}

\noindent In Table \ref{table:adaptive_refinement_nse}, we display the results from five DWR refinement loops with a mixed order error estimator for the Navier-Stokes 2D-3 problem. Although we do not employ a divergence-free projection, the effectivity indices hint at a good match between the true error and the error estimator. However, in the third loop a significant underestimation of the error in the mean drag can be observed which nevertheless does not persist 
after further refinement.

\begin{table}[H]
    \centering
    \resizebox{\columnwidth}{!}{%
    \begin{tabular}{|r|r|r|c|c|c|r|r|}
      \hline
     \#DoF(primal) & \#DoF(adjoint) & $M$ & $\eta_k$ & $\eta_h$ & $\eta$ & $J(\bm{U}) - J(\bm{U}_{\bm{kh}})$ & $\Ieff$ \\ \hline
          17,800 &  96,600 & 20 & $-9.3298 \cdot 10^{-6}$ & $3.8542 \cdot 10^{-1}$ & $3.8541 \cdot 10^{-1}$ & $5.5752 \cdot 10^{-1}$ & 0.69 \\ \hline
          63,454 &  350,922 & 36 & $-1.4015 \cdot 10^{-7}$ & $2.6505 \cdot 10^{-1}$ & $2.6505 \cdot 10^{-1}$ & $2.6642 \cdot 10^{-1}$ & 0.99 \\ \hline
          230,032 &  1,294,482 & 64 & $8.9182 \cdot 10^{-4}$ & $-1.2571 \cdot 10^{-2}$ & $1.1679 \cdot 10^{-2}$ & $1.2586 \cdot 10^{-1}$ & 0.09 \\ \hline
          828,744 &  4,706,883 & 113 & $-1.1615 \cdot 10^{-1}$ & $7.6888 \cdot 10^{-2}$ & $3.9265 \cdot 10^{-2}$ & $2.5449 \cdot 10^{-2}$ & 1.54 \\ \hline
          3,004,686 &  17,251,722 & 199 & $4.3194 \cdot 10^{-3}$ & $1.9094 \cdot 10^{-2}$ & $2.3414 \cdot 10^{-2}$ & $1.9674 \cdot 10^{-2}$ & 1.19 \\ \hline
    \end{tabular}
    }
    \caption{Adaptive refinement of mixed order on dynamic meshes for Navier-Stokes 2D-3 with divergence-free $L^2$ projection}
    \label{table:adaptive_refinement_nse_L2}
\end{table}

\noindent In Table \ref{table:adaptive_refinement_nse_L2}, we observe that the accuracy of the error estimator is improved compared to Table \ref{table:adaptive_refinement_nse} by employing a divergence-free $L^2$ projection. Additionally, this reduces the true error and thus yields a faster convergence towards the true solution.

\begin{table}[H]
    \centering
    \resizebox{\columnwidth}{!}{%
    \begin{tabular}{|r|r|r|c|c|c|r|r|}
      \hline
     \#DoF(primal) & \#DoF(adjoint) & $M$ & $\eta_k$ & $\eta_h$ & $\eta$ & $J(\bm{U}) - J(\bm{U}_{\bm{kh}})$ & $\Ieff$ \\ \hline
          17,800 &  96,600 & 20 & $-9.1676 \cdot 10^{-6}$ & $3.8548 \cdot 10^{-1}$ & $3.8547 \cdot 10^{-1}$ & $5.5768 \cdot 10^{-1}$ & 0.69 \\ \hline
          63,690 &  352,278 & 36 & $1.4710 \cdot 10^{-6}$ & $2.6493 \cdot 10^{-1}$ & $2.6493 \cdot 10^{-1}$ & $2.6667 \cdot 10^{-1}$ & 0.99 \\ \hline
          234,878 &  1,322,502 & 64 & $7.5795 \cdot 10^{-4}$ & $-4.4232 \cdot 10^{-3}$ & $3.6653 \cdot 10^{-3}$ & $1.2633 \cdot 10^{-1}$ & 0.03 \\ \hline
          834,710 &  4,741,881 & 113 & $-2.3546 \cdot 10^{-3}$ & $-7.6972 \cdot 10^{-2}$ & $7.9327 \cdot 10^{-2}$ & $1.9651 \cdot 10^{-2}$ & 4.04 \\ \hline
          3,044,708 &  17,485,449 & 199 & $3.0977 \cdot 10^{-3}$ & $9.0900 \cdot 10^{-3}$ & $1.2188 \cdot 10^{-2}$ & $6.5227 \cdot 10^{-3}$ & 1.87 \\ \hline
    \end{tabular}
    }
    \caption{Adaptive refinement of mixed order on dynamic meshes for Navier-Stokes 2D-3 with divergence-free $H^1_0$ projection}
    \label{table:adaptive_refinement_nse_H1}
\end{table}

\noindent In Table \ref{table:adaptive_refinement_nse_H1}, we have the same setup as in Table \ref{table:adaptive_refinement_nse_L2} but now use a divergence-free $H^1_0$ projection. In comparison to Table \ref{table:adaptive_refinement_nse_L2}, the effectivity indices are now slightly worse, but the true error is lower.

\begin{figure}[H]%
    \centering
    \includegraphics{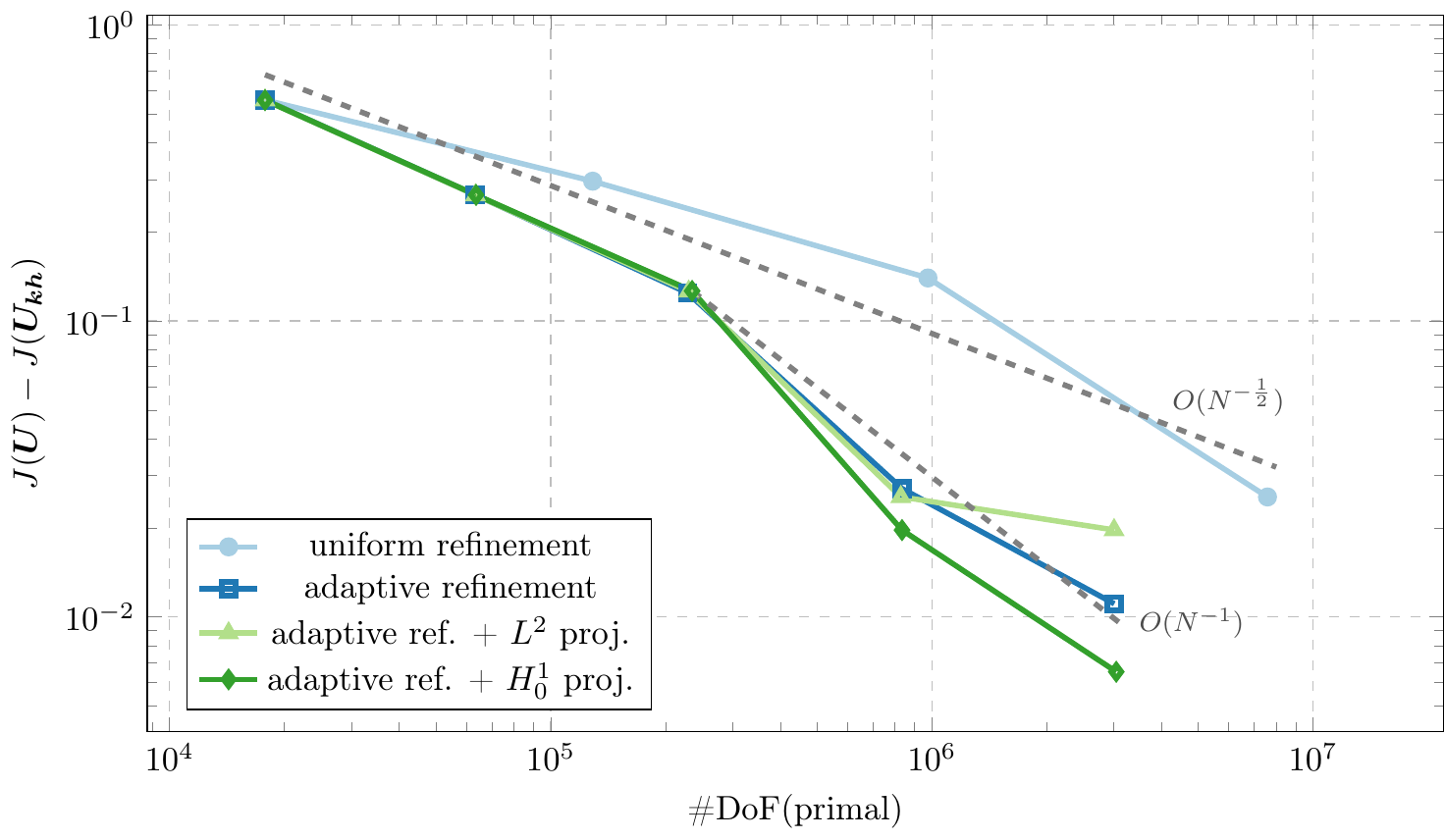}
    \caption{Comparison of the error in the mean drag for the Navier-Stokes problem for uniform and adaptive spatio-temporal refinement}%
    \label{fig:nse_comparison_error_uniform_and_adaptive_refinement}
\end{figure}

\noindent In Figure \ref{fig:nse_comparison_error_uniform_and_adaptive_refinement}, we visualize the true error from Tables \ref{table:uniform_refinement_nse}, \ref{table:adaptive_refinement_nse}, \ref{table:adaptive_refinement_nse_L2} and \ref{table:adaptive_refinement_nse_H1}. Analogous to Figure \ref{fig:stokes_comparison_error_uniform_and_adaptive_refinement}, we detect that the error of our adaptively refined solution converges faster than a uniformly refined solution.

\begin{figure}[H]%
    \centering
    \subfloat[Drag coefficient]{
        \includegraphics[clip,trim=3 0 4 0,     width=0.45\textwidth]{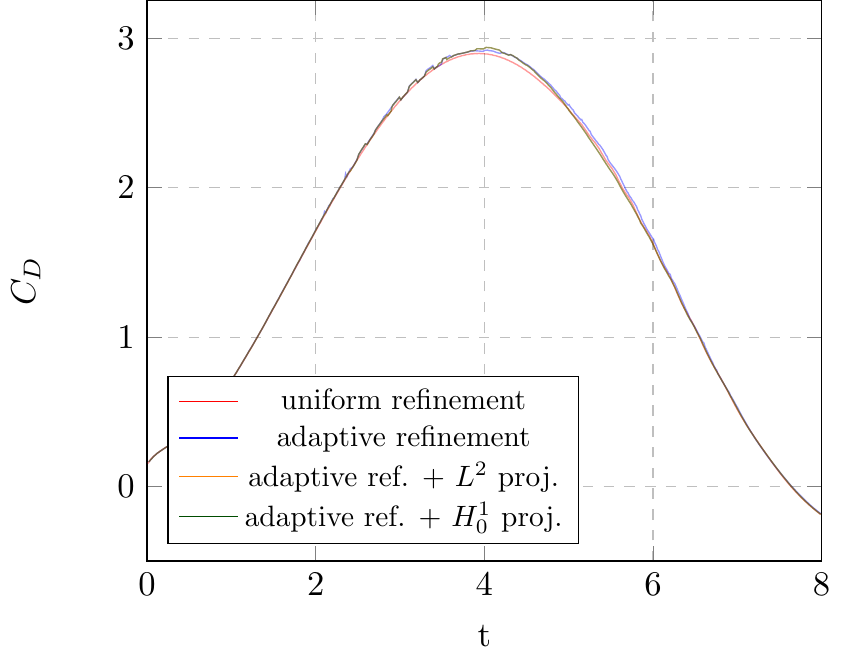}
    }\qquad 
    \subfloat[Lift coefficient]{
        \includegraphics[clip,trim=3 0 4 0,     width=0.45\textwidth]{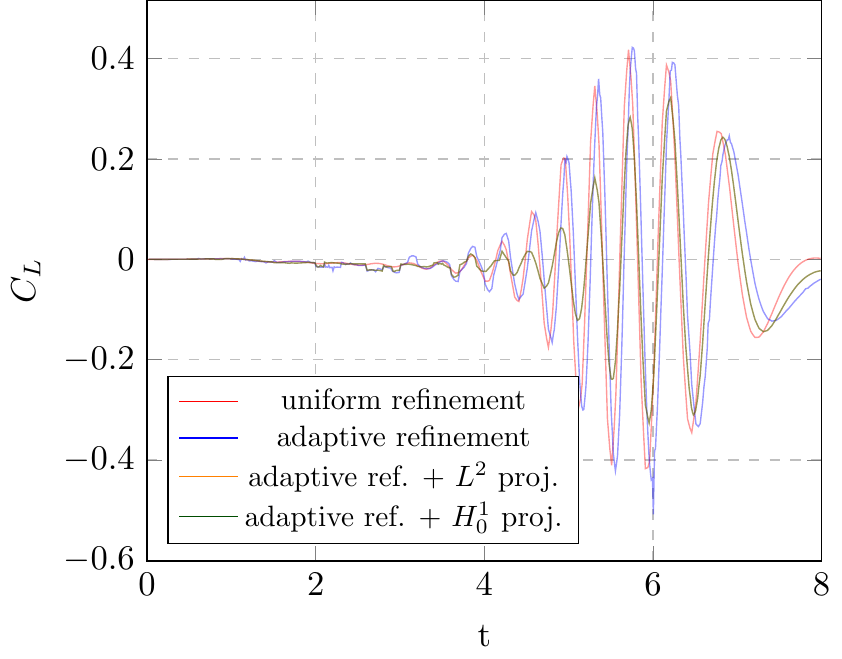}
    }
    \caption{Trajectory of the drag/lift coefficient for the Navier-Stokes problem for adaptive refinement}%
    \label{fig:nse_goal_functionals_adaptive}
\end{figure}

\noindent In Figure \ref{fig:nse_goal_functionals_adaptive}, we compare the trajectory of the drag and lift coefficient for the Navier-Stokes 2D-3 benchmark problem under different kinds of adaptive refinement. Like we observed in Figure \ref{fig:stokes_goal_functionals_adaptive}, the divergence free projections create less oscillatory drag and lift trajectories. This can be seen especially for the lift coefficient, where the projections help to remove the large spikes caused by using dynamic meshes. At the same time this causes a slight damping between $t = 4\ \text{s}$ and $t = 6\ \text{s}$. 

\begin{figure}[H]%
    \centering
    \includegraphics[scale=0.4]{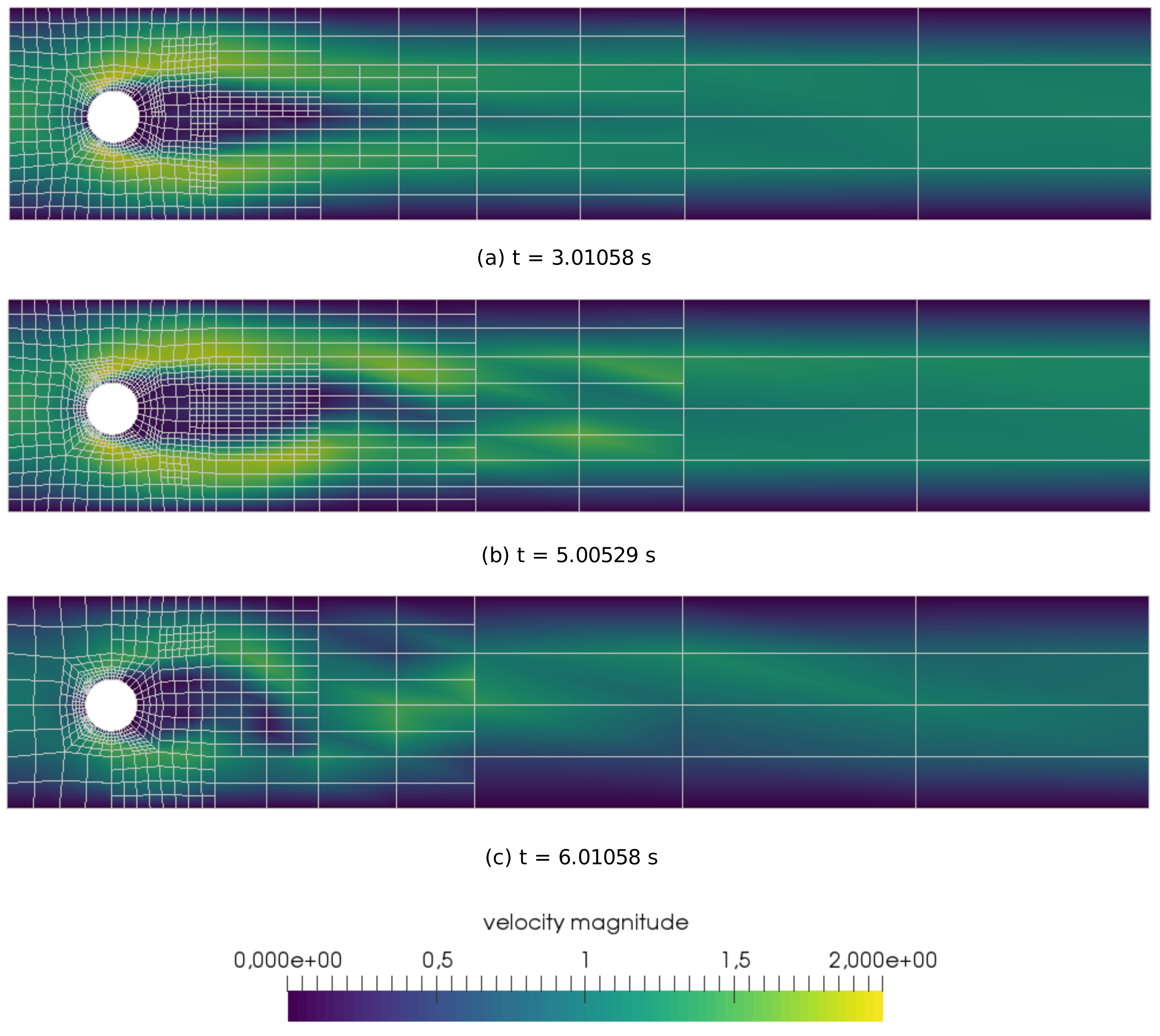}
    \caption{Snapshots of the velocity magnitude of the primal Navier-Stokes problem after 4 adaptive refinements with divergence-free $L^2$ projection}%
    \label{fig:nse_adaptive_velocity}
\end{figure}

\noindent In Figure \ref{fig:nse_adaptive_velocity}, we show some dynamic adaptive spatial meshes together with the snapshots of the velocity magnitude of the primal Navier-Stokes problem after 4 adaptive refinements with divergence-free $L^2$ projection. We observe that mainly the region around the cylinder is being refined, which is to be expected given that the drag coefficient is a boundary integral. Further away from the cylinder, e.g. close to the outflow boundary, the solution does not need to be resolved with the same accuracy. Moreover, since we are using dynamic meshes, the spatial mesh refinement can change over time as shown in Figure  \ref{fig:nse_adaptive_velocity}.

\begin{figure}[H]
    \centering
    \includegraphics{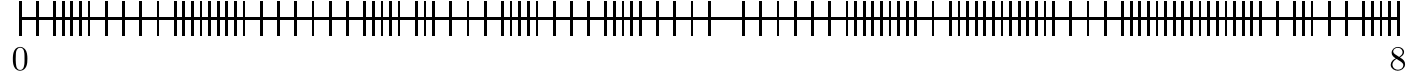}
    \caption{Temporal mesh in the final DWR loop for Navier-Stokes 2D-3 with divergence-free $L^2$ projection}
    \label{fig:nse_final_time_mesh}
\end{figure}

\noindent In Figure \ref{fig:nse_final_time_mesh}, we show the temporal mesh in the fifth DWR loop for the Navier-Stokes 2D-3 benchmark when using a divergence-free $L^2$ projection. The temporal mesh is also being refined adaptively and we can see that in some regions the diameter of the temporal elements is smaller, especially in the second half of the temporal domain.

\subsection{Navier-Stokes flow over backward facing step}

\begin{figure}[H]
    \begin{center}
    \includegraphics{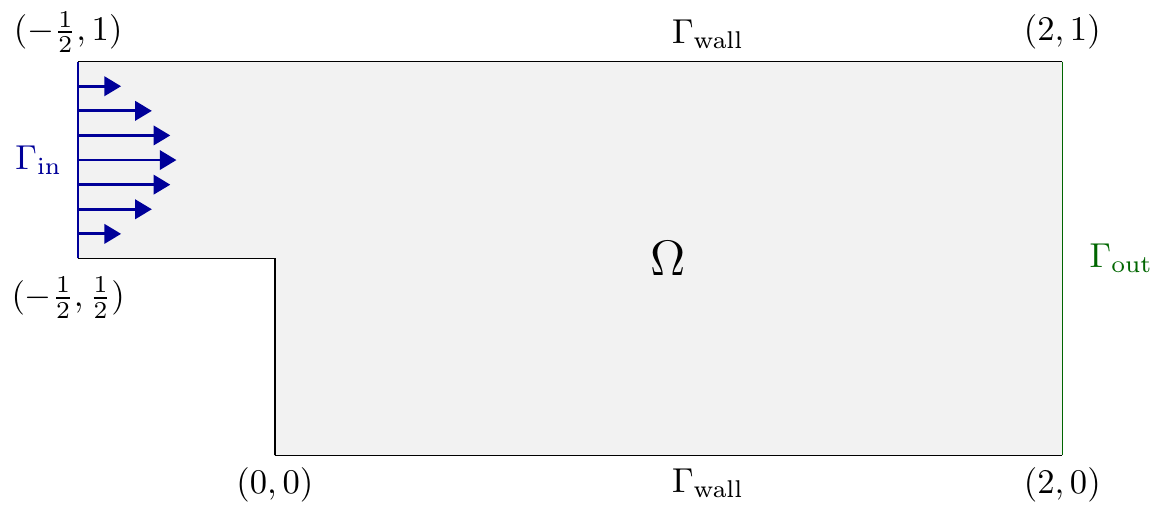}
    \caption{Domain of the 2D backward facing step benchmark}
    \label{fig:backward_step_domain}
    \end{center}
\end{figure}

In the third numerical experiment, we consider a version of the two-dimensional laminar backward-facing step benchmark \cite{Armaly1983, Biswas2004} with a time-dependent inflow condition. The domain is shown in Figure \ref{fig:backward_step_domain}.

\noindent The domain is defined as $\Omega := (-\frac{1}{2},0) \times (\frac{1}{2},1) \cup (0,2) \times (0,1)$.
The kinematic viscosity is $\nu = 1$ and we consider the time interval $I = (0,8)$. We prescribe homogeneous Neumann boundary conditions on the outflow boundary $\Gamma_{\text{out}} := \{2\} \times (0,1)$. On the inflow boundary $\Gamma_{\text{in}} := \{-\frac{1}{2}\} \times (\frac{1}{2},1)$ we enforce the time-dependent parabolic inflow profile
\begin{align}
	\bm{v}_D\left(t,-\frac{1}{2},y\right) = \begin{pmatrix}
		\sin\left(\frac{\pi t}{8}\right)80\left(-1+(3-2y)y\right) \\
		0
	\end{pmatrix}.
\end{align}
On all remaining boundaries homogeneous Dirichlet boundary conditions are being applied.
In contrast to the previous test cases, we now consider the mean squared vorticity
\begin{align}
    J_{\text{vorticity}}(\bm{U}) = \frac{1}{8}\int_0^8 \int_\Omega (\nabla \times \bm{v})^2\ \mathrm{d}\bm{x} \ \mathrm{d}t = \frac{1}{8}\int_0^8 \int_\Omega (\partial_1 v_2 - \partial_2 v_1)^2\ \mathrm{d}\bm{x}\ \mathrm{d}t
\end{align}
 as our goal functional. Note that this is a nonlinear goal functional and the right hand side of the dual problem is given by 
\begin{align}
    J^\prime_{\text{vorticity}, \bm{U}}(\bm{U})(\bm{\delta U}) = \frac{2}{8}\int_0^8 \int_\Omega (\nabla \times \bm{v}) \cdot (\nabla \times \bm{\delta v})\ \mathrm{d}\bm{x} \ \mathrm{d}t.
\end{align}
Using 80 initial temporal degrees of freedom and a spatial grid with 
1,439 primal degrees of freedom, we obtained the series of functional values from Table \ref{table:reference_bfs}.
We can see that this series converges from below and that the accuracy is
increasing with each refinement step.
\begin{table}[H]
    \centering
    \resizebox{0.7\columnwidth}{!}{%
    \begin{tabular}{|r|r|c|}
      \hline
     \#temporal DoF(primal) & \#spatial DoF(primal) & $J(\bm{U}_{\bm{kh}})$
     \\ \hline
          80
          & 1,439
          & 506.20900857749609
          \\ \hline 
          160
          & 5,467
          & 510.27887058685025
          \\ \hline
          320 
          & 21,299
          & 512.24259375586087
          \\ \hline
          640
          & 84,067
          & 513.18333123881553
          \\ \hline
          1,280 
          & 334,019
          & 513.63122326525161
          \\ \hline
          2,560
          & 1,331,587 
          & 513.84343972465513
          \\ \hline
    \end{tabular}
    }
    \caption{Mean vorticity value under uniform spatio-temporal refinement for the backward-facing step test case}
    \label{table:reference_bfs}
\end{table}

\noindent For the sake of brevity we will only discuss adaptive results in detail.
For this test case we used the fixed rate marking strategy both in space and time
with $\varepsilon = 30$.

\begin{table}[H]
    \centering
    \resizebox{0.85\columnwidth}{!}{%
    \begin{tabular}{|c||c|c||c|c||c|c|}
      \hline 
     DWR loop & $\text{coarse} \Ieff$ & $\text{fine} \Ieff$ 
     &  $\text{coarse} \Ieff^{L_2}$ &  $\text{fine} \Ieff^{L_2}$ 
     & $ \text{coarse} \Ieff^{H^1_0}$ & $ \text{fine} \Ieff^{H^1_0}$
     \\ \hline
        1 & 1.06 & 1.03 & 1.05 & 1.03 &  1.06 & 1.03 \\ \hline
        2 & 1.15 & 1.08 & 1.15 & 1.08 &  1.10 & 1.04 \\ \hline
        3 & 1.34 & 1.17 & 1.34 & 1.17 &  1.83 & 1.53 \\ \hline
        4 & 2.02 & 1.37 & 2.21 & 1.47 &  -2.51 & -3.97 \\ \hline
        5 & -1.11 & 2.06 & -1.53 & 2.47 &  -6.46 & -0.70 \\ \hline
    \end{tabular}
    }
    \caption{Signed effectivity indices for the different projections with
    a once coarsened spatio-temporal mesh (coarse) compared to the 
    finest reference mesh (fine). Positive $\Ieff$ means the computed functional value
    is below the reference value.}
    \label{table:effectivity_bfs}
\end{table}
\noindent Table \ref{table:effectivity_bfs} shows the resulting effectivity indices 
for the different projection strategies (no projection, $L^2$ projection, $H^1_0$ projection). 
To investigate the effect of the reference value on the approximated real error,
we compare the results with the last (fine) and second to last (coarse) 
mean vorticity as reference values $J_{\text{vorticity}}(\bm{U})$. 
As long as the vorticity on the adaptive mesh is below the reference value
we can see that we get a better effectivity with the finer reference value
suggesting that our estimator is a good predictor of the true error. Additionally, we see a bad effectivity index for the loop in which
we have the same $h_{\min}$ as in the reference grid. 
However, this is not surprising as it is generally not advisable to analyze 
adaptive results on the reference mesh level.

\begin{figure}[H]%
    \centering
    \includegraphics{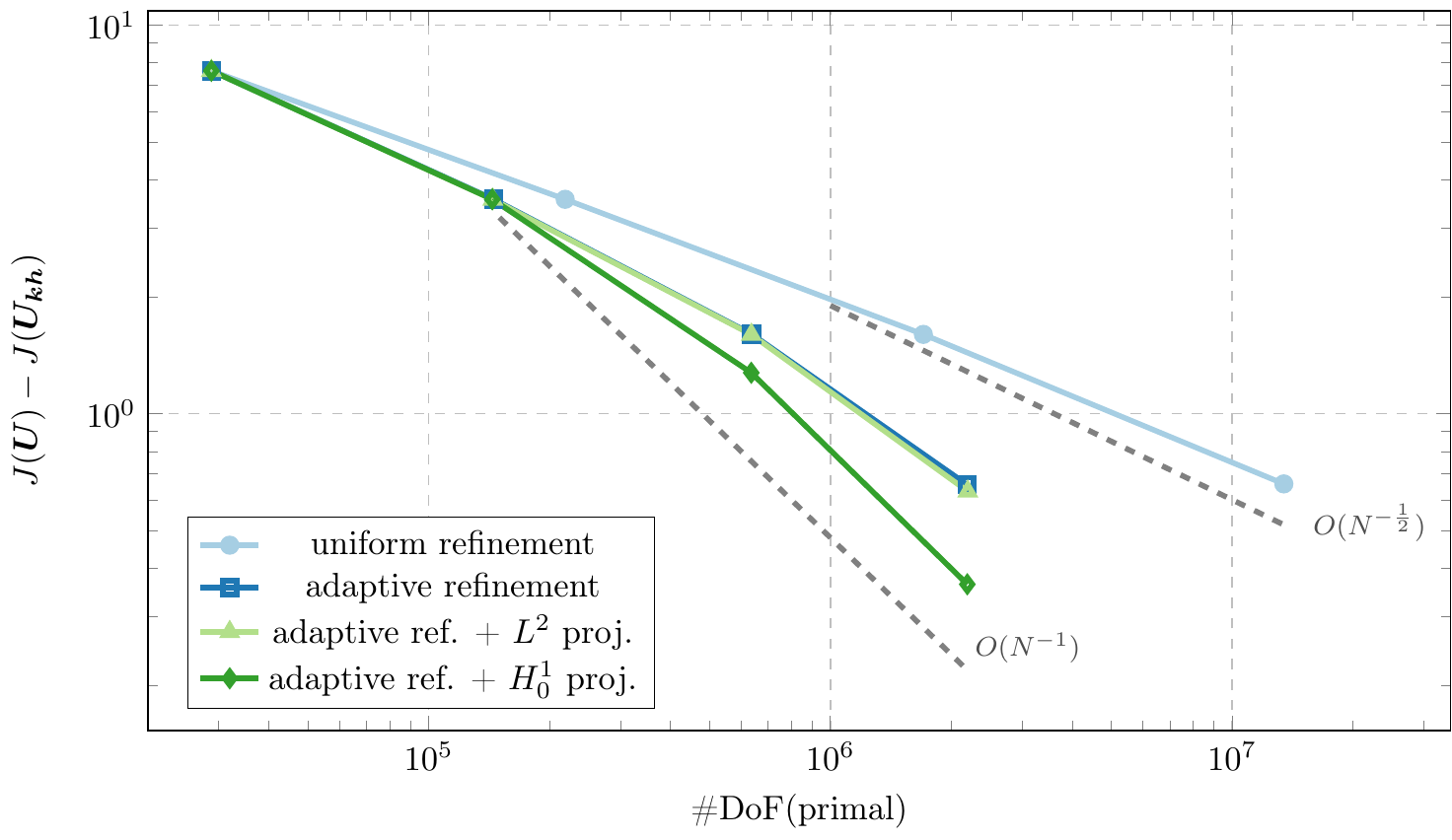}
    \caption{Comparison of the true error in the mean squared vorticity for the backward facing step for uniform and adaptive spatio-temporal refinement}%
    \label{fig:backward_step_comparison_error_uniform_and_adaptive_refinement}
\end{figure}

\noindent In Figure \ref{fig:backward_step_comparison_error_uniform_and_adaptive_refinement}, we visualize the true error in the mean squared vorticity for the backward facing step for uniform and adaptive spatio-temporal refinement. Analogous to Figure \ref{fig:stokes_comparison_error_uniform_and_adaptive_refinement} and Figure \ref{fig:nse_comparison_error_uniform_and_adaptive_refinement} for the (Navier-)Stokes 2D-3 benchmark, we detect that the error of our adaptively refined solution converges with a higher rate than a uniformly refined solution. In particular, when comparing the fourth refinement cycle, we observe that adaptive refinement uses only $\frac{1}{6}$ of the number of space-time degrees of freedom as uniform refinement to reach a tolerance of $0.7$ in the goal functional.

\begin{figure}[H]%
    \centering
    \includegraphics[width=0.7\textwidth]{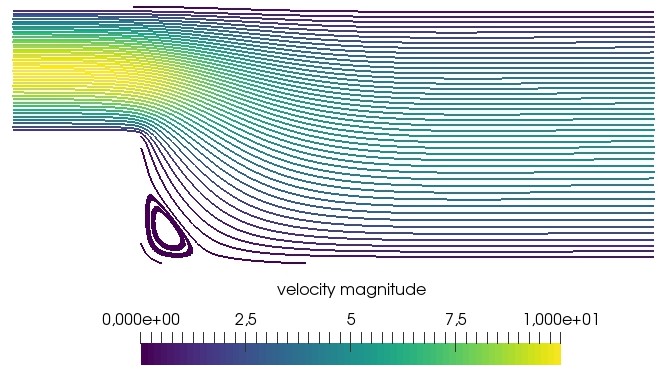}
    \caption{Snapshot of the velocity magnitude and streamlines of the primal backward-facing step problem at $t\approx4$}%
    \label{fig:bfs_velocity_mag}
\end{figure}
\noindent In Figure \ref{fig:bfs_velocity_mag}, we show the snapshot of the velocity magnitude of the primal backward-facing step problem at time $ t \approx 4$ and its streamlines. Looking at the streamlines, we observe a corner vortex in the lower left corner after the step.

\begin{figure}[H]%
    \centering
    \includegraphics[width=0.7\textwidth]{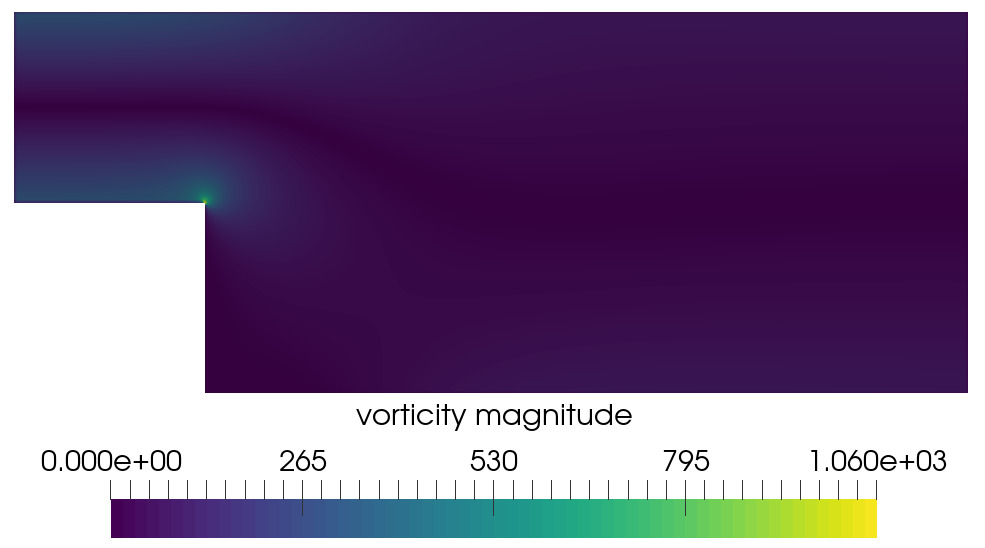}
    \caption{Snapshot of the vorticity magnitude of the primal backward-facing step problem at $t\approx4$}%
    \label{fig:bfs_vorticity_mag}
\end{figure}

\noindent In Figure \ref{fig:bfs_vorticity_mag}, we plot the snapshot of the vorticity magnitude of the primal backward-facing step problem at time $ t \approx 4$. It can be seen that the vorticity magnitude is largest at the reentrant corner $\bm{x} = (0,\frac{1}{2})$.

\begin{figure}[H]%
    \centering
    \includegraphics[width=0.7\textwidth]{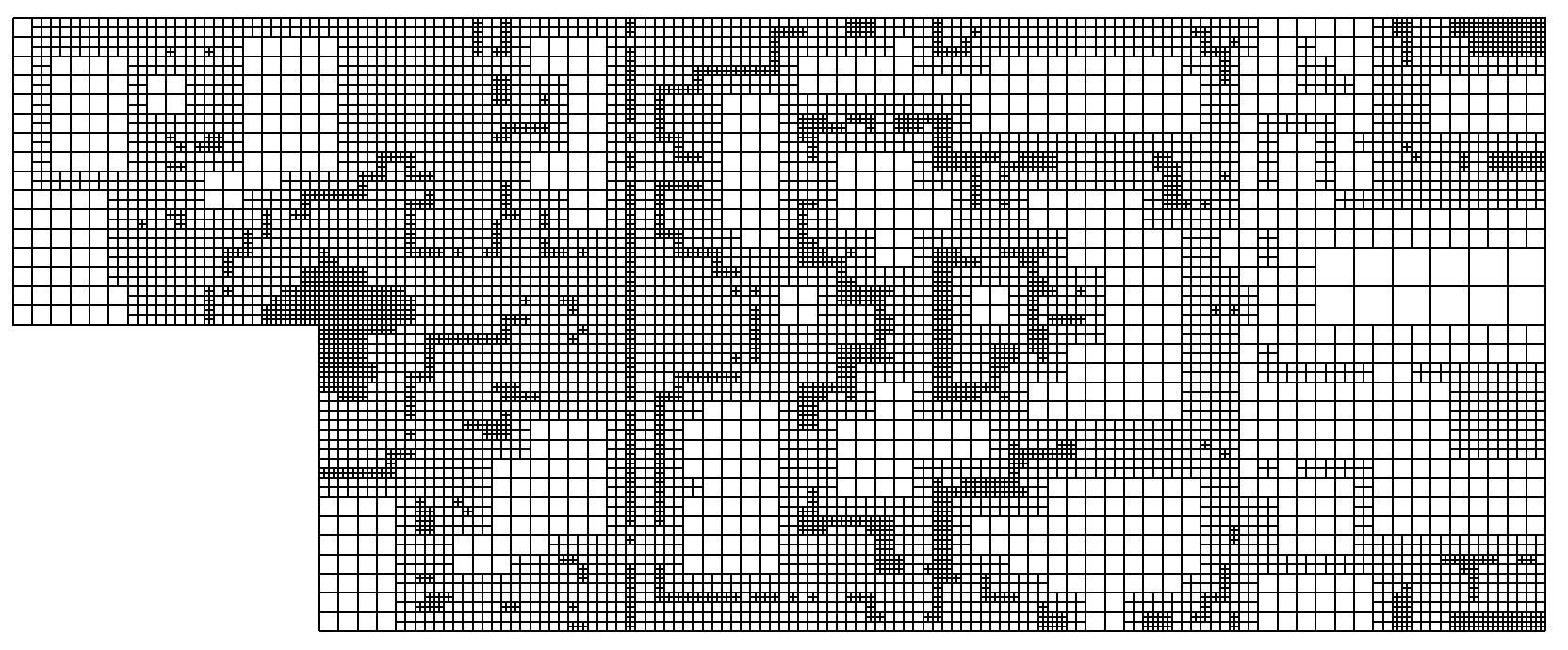}
    \caption{Spatial grid at $t\approx4$ for the backward-facing step problem after 4 adaptive refinements with divergence-free $L^2$ projection}%
    \label{fig:bfs_grid}
\end{figure}
\noindent In Figure \ref{fig:bfs_grid}, we display a 4 times adaptively refined grid with divergence-free $L^2$ projection. On the one hand, we observe that the grid is mostly being refined close to the reentrant corner where the vorticity magnitude is at its maximum. On the other hand, closer to the inflow or outflow boundaries the mesh is much coarser.

%%%%%%%%%%%%%%%%%%%%%%%%%%%%%%%%%%%%%%%%%%%%%%%%%%%%%%%%%%%%
%%                   CONCLUSION                           %%
%%%%%%%%%%%%%%%%%%%%%%%%%%%%%%%%%%%%%%%%%%%%%%%%%%%%%%%%%%%%
\section{Conclusions and outlook}
\label{sec_conclusions}
In this work, we formulated a space-time PU-DWR approach for
goal-oriented error control and local mesh adaptivity 
for the incompressible Navier-Stokes equations. 
Therein, our focus was temporal and spatial error control, which has been computationally 
analyzed for different benchmark problems such as the Sch\"afer-Turek 1996 benchmark.
In terms of convergence as well as evaluation of effectivity indices we observed excellent 
performances for different goal functionals.
In ongoing work, an interesting aspect would be to implement Linke's pressure-robust 
discretizations rather than the global corrections proposed by Besier and Wollner.

%%%%%%%%%%%%%%%%%%%%%%%%%%%%%%%%%%%%%%%%%%%%%%%%%%%%%%%%%%%%
%%                  ACKNOWLEDGEMENTS                      %%
%%%%%%%%%%%%%%%%%%%%%%%%%%%%%%%%%%%%%%%%%%%%%%%%%%%%%%%%%%%%
\section*{Acknowledgements}
The first and fourth author acknowledge the funding of the German Research Foundation (DFG) within the framework of the International Research Training Group on  Computational Mechanics Techniques in High Dimensions GRK 2657 under Grant Number 433082294. In addition, we thank Hendrik Fischer, Marius Paul Bruchhäuser, Bernhard Endtmayer and Ludovic Chamoin for fruitful discussions and comments.
%%%%%%%%%%%%%%%%%%%%%%%%%%%%%%%%%%%%%%%%%%%%%%%%%%%%%%%%%%%%

\bibliographystyle{abbrv}
%\bibliography{lit}

\end{document}